 \documentstyle[12pt]{article}
    \title{{\bf  Virasoro vertex operator algebras, the (nonmeromorphic)
operator product expansion and the tensor
product theory}}
    \author{Yi-Zhi Huang\thanks{Supported in part by NSF grant
DMS-9301020 and by DIMACS, an
NSF Science and Technology Center funded under contract STC-88-09648.}}
    \date{}
    
\begin{document}
    \bibliographystyle{alpha}
    \maketitle

    \input amssym.def
    \input amssym
    \newtheorem{rema}{Remark}[section]
    \newtheorem{propo}[rema]{Proposition}
    \newtheorem{theo}[rema]{Theorem}
   \newtheorem{defi}[rema]{Definition}
    \newtheorem{lemma}[rema]{Lemma}
    \newtheorem{corol}[rema]{Corollary}
     \newtheorem{exam}[rema]{Example}
	\newcommand{\nno}{\nonumber}
	\newcommand{\lbar}{\bigg\vert}
	\newcommand{\p}{\partial}
	\newcommand{\dps}{\displaystyle}
	\newcommand{\bra}{\langle}
	\newcommand{\ket}{\rangle}
 \newcommand{\res}{\mbox{\rm Res}}
 \newcommand{\pf}{{\it Proof}\hspace{2ex}}
 \newcommand{\epf}{\hspace{2em}$\Box$}
 \newcommand{\epfv}{\hspace{1em}$\Box$\vspace{1em}}
\newcommand{\nord}{\mbox{\scriptsize ${\circ\atop\circ}$}}
\newcommand{\wt}{\mbox{\rm wt}\ }

 \makeatletter
\newlength{\@pxlwd} \newlength{\@rulewd} \newlength{\@pxlht}
\catcode`.=\active \catcode`B=\active \catcode`:=\active \catcode`|=\active
\def\sprite#1(#2,#3)[#4,#5]{
   \edef\@sprbox{\expandafter\@cdr\string#1\@nil @box}
   \expandafter\newsavebox\csname\@sprbox\endcsname
   \edef#1{\expandafter\usebox\csname\@sprbox\endcsname}
   \expandafter\setbox\csname\@sprbox\endcsname =\hbox\bgroup
   \vbox\bgroup
  \catcode`.=\active\catcode`B=\active\catcode`:=\active\catcode`|=\active
      \@pxlwd=#4 \divide\@pxlwd by #3 \@rulewd=\@pxlwd
      \@pxlht=#5 \divide\@pxlht by #2
      \def .{\hskip \@pxlwd \ignorespaces}
      \def B{\@ifnextchar B{\advance\@rulewd by \@pxlwd}{\vrule
         height \@pxlht width \@rulewd depth 0 pt \@rulewd=\@pxlwd}}
      \def :{\hbox\bgroup\vrule height \@pxlht width 0pt depth
0pt\ignorespaces}
      \def |{\vrule height \@pxlht width 0pt depth 0pt\egroup
         \prevdepth= -1000 pt}
   }
\def\endsprite{\egroup\egroup}
\catcode`.=12 \catcode`B=11 \catcode`:=12 \catcode`|=12\relax
\makeatother

\def\hboxtr{\FormOfHboxtr} 
\sprite{\FormOfHboxtr}(25,25)[0.5 em, 1.2 ex] 

:BBBBBBBBBBBBBBBBBBBBBBBBB |
:BB......................B |
:B.B.....................B |
:B..B....................B |
:B...B...................B |
:B....B..................B |
:B.....B.................B |
:B......B................B |
:B.......B...............B |
:B........B..............B |
:B.........B.............B |
:B..........B............B |
:B...........B...........B |
:B............B..........B |
:B.............B.........B |
:B..............B........B |
:B...............B.......B |
:B................B......B |
:B.................B.....B |
:B..................B....B |
:B...................B...B |
:B....................B..B |
:B.....................B.B |
:B......................BB |
:BBBBBBBBBBBBBBBBBBBBBBBBB |

\endsprite

\def\shboxtr{\FormOfShboxtr} 
\sprite{\FormOfShboxtr}(25,25)[0.3 em, 0.72 ex] 

:BBBBBBBBBBBBBBBBBBBBBBBBB |
:BB......................B |
:B.B.....................B |
:B..B....................B |
:B...B...................B |
:B....B..................B |
:B.....B.................B |
:B......B................B |
:B.......B...............B |
:B........B..............B |
:B.........B.............B |
:B..........B............B |
:B...........B...........B |
:B............B..........B |
:B.............B.........B |
:B..............B........B |
:B...............B.......B |
:B................B......B |
:B.................B.....B |
:B..................B....B |
:B...................B...B |
:B....................B..B |
:B.....................B.B |
:B......................BB |
:BBBBBBBBBBBBBBBBBBBBBBBBB |

\endsprite

\begin{abstract}
In \cite{HL1}--\cite{HL5} and \cite{H1}, a theory of tensor products of
modules for a vertex operator algebra is being developed. To use this theory,
one first has to verify that the vertex operator algebra satisfies
 certain conditions. We show in the present paper that
for any vertex operator algebra containing a vertex operator
subalgebra isomorphic to a tensor product algebra of  minimal
Virasoro vertex operator algebras (vertex
operator algebras associated to minimal models),
the tensor product theory can be applied. In particular, intertwining
operators for such a vertex operator algebra satisfy the
(nonmeromorphic) commutativity
(locality) and the (nonmeromorphic) associativity (operator product expansion).
Combined  with a result announced in
\cite{HL4}, the results of the present paper also show that the category
of modules for such a vertex operator algebra has a natural structure of
a braided tensor category. In particular, for
 any pair $p, q$ of relatively prime positive integers larger
than $1$, the category of
minimal modules of central charge $1-6\frac{(p-q)^{2}}{pq}$
for the Virasoro algebra has a natural structure of a braided tensor category.
\end{abstract}

\renewcommand{\theequation}{\thesection.\arabic{equation}}
\renewcommand{\therema}{\thesection.\arabic{rema}}
\setcounter{equation}{0}
\setcounter{rema}{0}
\setcounter{section}{-1}

\section{Introduction}

Vertex operator algebras were introduced by Borcherds \cite{B} and, in
a modified form, by Frenkel, Lepowsky and Meurman \cite{FLM}.  They
are essentially the same as chiral algebras in physics (see, for
example, \cite{MS}).  One crucial property of vertex operator algebras
and modules for them is the meromorphicity.  It is the meromorphicity
that makes it possible to express the axioms of vertex operator
algebras using formal series or components of vertex operators. These
axioms expressed using formal series or components of vertex operators
allow us to construct many examples of vertex operator algebras using
algebraic methods. But in conformal field theory and in applications
of conformal field theory, the more interesting objects to study are
intertwining operators (or chiral vertex operators). In general,
intertwining operators do not satisfy the meromorphicity condition since
by definition they involve noninteger powers of a variable. The most
fundamental assumption on intertwining operators is the so called
``operator product expansion'' in physics. In the terminology of the
theory of vertex operator algebras, the operator product expansion
of chiral vertex operators is
exactly the associativity of intertwining operators.  Many beautiful
results in conformal field theory are obtained based on this
assumption. One example is the braided tensor category structures
constructed {}from conformal field theories \cite{MS}. In \cite{MS},
Moore and Seiberg obtained the braiding and fusing matrices (the
commutativity and the associativity isomorphisms) using the assumed
(nonmeromorphic) operator product expansion of chiral vertex operators
(intertwining
operators) or an equivalent geometric axiom of conformal field
theories.

Another example is the study of orbifold conformal field theories \cite{DHVW1}
\cite{DHVW2}
\cite{DFMS} \cite{DVVV}.
Orbifold conformal field theories (or simply orbifold theories)
are a rich source of conformal
field theories. The moonshine module constructed by Frenkel, Lepowsky and
Meurman \cite{FLM1} \cite{FLM} is historically the first example
of orbifold theories. Dolan, Goddard and Montague gave another proof
that the moonshine module is a ``meromorphic conformal field theory''
(vertex operator algebra) using techniques developed in string theory
\cite{DGM}. Their proof works for
a class of ${\Bbb Z}_{2}$-orbifold theories  and thus allows
them to obtain a better understanding about Frenkel-Lepowsky-Meurman's triality
\cite{DGM2} \cite{DGM3}. Another orbifold construction---a
$\Bbb{Z}_{3}$-orbifold construction---of the moonshine module
was given by Dong and Mason \cite{DM1}. In \cite{Tu1} and
\cite{Tu2}, Tuite showed that the Monstrous moonshine
conjectured by Conway and Norton \cite{CN} and proved by Borcherds \cite{B2}
for Frenkel-Lepowsky-Meurman's moonshine module
can be understood using the uniqueness,  conjectured
by Frenkel, Lepowsky and Meurman \cite{FLM}, of the moonshine module
and some conjectures on the
orbifold theories constructed {}from the moonshine module. Orbifold
 theories were also used to construct the mirror of a Calabi-Yau manifold
\cite{GP}. In all of the studies of orbifold theories mentioned above,
the results are either only about meromorphic parts or
based on the basic assumption that the nonmeromorphic operator product
expansion holds.

To establish conformal field theory as a solid mathematical theory and
to solve concrete mathematical problems using conformal field
theories, we need to prove, not assume, the associativity (operator
product expansion) of intertwining operators (chiral vertex
operators). The main difficulty---but also the source of the richness
of conformal field theories---is that products and iterates of
intertwining operators are in general not meromorphic.  One way to
do this is to construct the intertwining operators directly {}from
relatively elementary mathematical data and to then show that they
satisfy associativity.  This method seems to be very difficult and is
also not practical since one has to construct intertwining operators
and to prove the associativity case by case. Another method, the
method that we use in the present and related papers, is to study
intertwining operators as objects in the representation theory of
vertex operator algebras.

In the representation theory of Lie algebras, one of the most important
operations is the tensor product operation for modules for Lie algebras.
In the representation theory of vertex operator algebras, a tensor product
theory for modules for a vertex operator algebra is being developed in
\cite{HL1}--\cite{HL5} and \cite{H1}. But to use this theory,
one first has to verify that the vertex operator algebra satisfies
certain conditions. The purpose of the present paper is to show that
for any vertex operator algebra containing a vertex operator
subalgebra isomorphic to a tensor product algebra of minimal Virasoro
vertex operator algebras (vertex operator algebras associated minimal
models), the tensor product theory can be applied. In particular,
intertwining operators for such a vertex operator algebra satisfy the
(nonmeromorphic) commutativity (locality) and the (nonmeromorphic)
associativity (operator product expansion).  Recall that a vertex operator
subalgebra of a vertex operator algebra is a subspace which contains the vacuum
vector and the Virasoro element and is closed under the vertex operator
map \cite{FHL}. The condition that a vertex operator algebra contains a
vertex operator
subalgebra isomorphic to a tensor product algebra of minimal Virasoro
vertex operator algebras gives enough information we need
about the vertex operator
algebra. For example, if the vertex operator algebra contains
a trivial one-dimensional vertex operator
subalgebra isomorphic to the null tensor product of minimal Virasoro
vertex operator algebras, then
the vertex operator algebra is a finite-dimensional
commutative associative algebra.

Combined with a result
announced in
\cite{HL4}, the results of the present paper also show that the category
of modules for such a vertex operator algebra has a natural structure of
a braided tensor category. In particular, for
 any pair $p, q$ of relatively prime positive integers larger
than $1$, the category of
minimal modules of central charge $1-6\frac{(p-q)^{2}}{pq}$
for the Virasoro algebra has a natural structure of a braided tensor category.
Since in general
the vertex operator algebras studied in this paper are not rational
 (see Example \ref{irra}), we also relax the conditions to
use the tensor product theory given in \cite{HL2} \cite{HL4} \cite{HL5} and
\cite{H1}, especially the rationality of the vertex operator algebra, in
this paper.
For the precise statements of the results in the present paper, see Section 3.

The results of the present paper have many applications.
One immediate consequence is that for many vertex operator algebras
associated to ${\cal W}$-algebras, the
tensor product theory can be applied.
The results of the present paper
 have been used in the proof by the author that there is a natural structure
of an abelian intertwining algebra (in the sense of Dong and Lepowsky
\cite{DL}) on the direct sum of the untwisted
vertex operator
algebra constructed {}from the Leech lattice and its (unique)
irreducible twisted module (see \cite{H2}).
A special case of this abelian intertwining algebra
structure gives a new and conceptual proof that the moonshine module is
a vertex operator algebra. This abelian intertwining algebra also contains as
substructures several other interesting structures, including a
twisted module for the
moonshine module, the superconformal structure of Dixon, Ginsparg and Harvey
\cite{DGH}, a vertex operator superalgebra and a twisted module for it.
We believe that the results of the
present paper can also be used to give a new and conceptual proof of
the result of
Dolan, Goddard and Montague  \cite{DGM}.
The results of the present
paper can also be reformulated in terms of (partial) operads, and genus-zero
conformal field theories in the sense of Segal \cite{S}
can be constructed using this
reformulation. For vertex operator algebras containing as vertex operator
subalgebras tensor
product algebras of vertex operator algebras associated to affine Lie algebras,
we have similar results. All of these will be discussed in separate
papers.

This paper is organized as follows:
In Section 1, we  review briefly
the tensor product theory for
modules for a vertex operator algebra developed in
\cite{HL1}--\cite{HL5} and \cite{H1}. In Section 2, we
review the results on the minimal Virasoro vertex operator algebras
and their
modules obtained in \cite{FZ}, \cite{DMZ} and \cite{W}.  We state the
results of the present paper, including Theorem \ref{commu}, Theorem
\ref{subvtc}, Theorem \ref{cpqm}, Proposition \ref{cpqm1}, Theorem
\ref{main}, Corollary \ref{cpqmbtc} and Theorem \ref{virbtc}, in
Section 3.  In Section 3, we also give an example showing that for any
rational vertex operator algebra, there exists an irrational vertex
operator algebra containing this rational vertex operator algebra as a
vertex operator subalgebra.  Theorem \ref{commu}, Theorem
\ref{subvtc}, Theorem \ref{cpqm} and Proposition \ref{cpqm1} are
proved in Section 4, Section 5, Section 6 and Section 7, respectively.

{\bf Acknowledgement} I would like to thank James Lepowsky for
comments.

\renewcommand{\theequation}{\thesection.\arabic{equation}}
\renewcommand{\therema}{\thesection.\arabic{rema}}
\setcounter{equation}{0}
\setcounter{rema}{0}

\section{The tensor product theory}

We review briefly in this section the tensor product theory for
modules for a vertex operator algebra being developed in
\cite{HL1}--\cite{HL5} and
\cite{H1}.

Fix $z\in \Bbb{C}^{\times}$ and let $(W_{1}, Y_{1})$, $(W_{2},
Y_{2})$ and $(W_{3}, Y_{3})$ be $V$-modules.  A {\it
$P(z)$-intertwining map of type ${W_{3}}\choose {W_{1}W_{2}}$} is a
linear map $F: W_{1}\otimes W_{2} \to
\overline{W}_{3}$ satisfying the condition
\begin{eqnarray*}
\lefteqn{x_{0}^{-1}\delta\left(\frac{ x_{1}-z}{x_{0}}\right)
Y_{3}(v, x_{1})F(w_{(1)}\otimes w_{(2)})=}\nonumber\\
&&=z^{-1}\delta\left(\frac{x_{1}-x_{0}}{z}\right)
F(Y_{1}(v, x_{0})w_{(1)}\otimes w_{(2)})\nonumber\\
&&\hspace{2em}+x_{0}^{-1}\delta\left(\frac{z-x_{1}}{-x_{0}}\right)
F(w_{(1)}\otimes Y_{2}(v, x_{1})w_{(2)})
\end{eqnarray*}
for $v\in V$, $w_{(1)}\in W_{1}$, $w_{(2)}\in W_{2}$.  The vector
space of $P(z)$-intertwining maps of type ${W_{3}}\choose
{W_{1}W_{2}}$ is denoted by ${\cal M}[P(z)]^{W_{3}}_{W_{1}W_{2}}$.  A
{\it $P(z)$-product of $W_{1}$ and $W_{2}$} is a $V$-module $(W_{3},
Y_{3})$ equipped with a $P(z)$-intertwining map $F$ of type
${W_{3}}\choose {W_{1}W_{2}}$ and is denoted by $(W_{3}, Y_{3}; F)$
(or simply by $(W_{3}, F)$).  Let $(W_{4}, Y_{4}; G)$ be another
$P(z)$-product of $W_{1}$ and $W_{2}$.  A {\it morphism} {}from
$(W_{3}, Y_{3}; F)$ to $(W_{4}, Y_{4}; G)$ is a module map $\eta$
{}from $W_{3}$ to $W_{4}$ such that $G=\overline{\eta}\circ F$, where
$\overline{\eta}$ is the map {}from $\overline{W}_{3}$ to
$\overline{W}_{4}$ uniquely extending $\eta$.

A {\it $P(z)$-tensor product of $W_{1}$ and $W_{2}$} is a
$P(z)$-product $$(W_{1}\boxtimes_{P(z)} W_{2}, Y_{P(z)};
\boxtimes_{P(z)})$$
such that for any $P(z)$-product
$(W_{3}, Y_{3}; F)$, there is a unique morphism {}from
$$(W_{1}\boxtimes_{P(z)} W_{2}, Y_{P(z)};
\boxtimes_{P(z)})$$ to $(W_{3}, Y_{3}; F)$.
The $V$-module $(W_{1}\boxtimes_{P(z)} W_{2},  Y_{P(z)})$ is
called a {\it $P(z)$-tensor product module} of $W_{1}$ and $W_{2}$.
A $P(z)$-tensor product is unique up to isomorphism.

To construct a $P(z)$-tensor product of $W_{1}$ and $W_{2}$,
we define an action of $$V \otimes \iota_{+}\Bbb{C}[t,t^{- 1},
(z^{-1}-t)^{-1}]$$ on $(W_1 \otimes W_2)^*$ (where
$\iota_{+}$ denotes the operation of expansion of a
rational function of $t$ in the direction of positive powers of $t$),
that is, a linear map
$$
\tau_{P(z)}: V\otimes \iota_{+}\Bbb{C}[t, t^{-1},
(z^{-1}-t)^{-1}]\to
\mbox{\rm End}\;(W_{1}\otimes W_{2})^{*},
$$
by
\begin{eqnarray*}
\lefteqn{\left(\tau_{P(z)}
\left(x_{0}^{-1}\delta\left(\frac{x^{-1}_{1}-z}{x_{0}}\right)
Y_{t}(v, x_{1})\right)\lambda\right)(w_{(1)}\otimes
w_{(2)})}\nonumber\\
&&=z^{-1}\delta\left(\frac{x^{-1}_{1}-x_{0}}{z}\right)
\lambda(Y_{1}(e^{x_{1}L(1)}(-x_{1}^{-2})^{L(0)}v, x_{0})w_{(1)}
\otimes w_{(2)})
\nonumber\\
&&\hspace{2em}+x^{-1}_{0}\delta\left(\frac{z-x^{-1}_{1}}{-x_{0}}
\right)
\lambda(w_{(1)}\otimes Y_{2}^{*}(v, x_{1})w_{(2)})
\end{eqnarray*}
for $v\in V$, $\lambda\in (W_{1}\otimes W_{2})^{*}$, $w_{(1)}\in
W_{1}$, $w_{(2)}\in W_{2}$, where
$$
Y_{t}(v, x)=v\otimes x^{-1}\delta\left(\frac{t}{x}\right).
$$
There is an obvious action of $$V \otimes \iota_{+}\Bbb{C}[t,t^{- 1},
(z^{-1}-t)^{-1}]$$ on any $V$-module. We have:

\begin{propo}\label{intwact}
Under the natural isomorphism
$$
\mbox{\rm Hom}(W'_{3}, (W_{1}\otimes
W_{2})^{*})\stackrel{\sim}{\to}\mbox{\rm Hom}(W_{1}\otimes W_{2},
\overline{W}_{3}),
$$
the maps in $\mbox{\rm Hom}(W'_{3}, (W_{1}\otimes
W_{2})^{*})$ intertwining the two actions of
$$V \otimes \iota_{+}\Bbb{C}[t,t^{-
1},(z^{-1}-t)^{-1}]$$ on $W'_{3}$ and
$(W_{1}\otimes W_{2})^{*}$ correspond exactly to the
$P(z)$-intertwining maps of type ${W_{3}}\choose {W_{1}W_{2}}$.
\end{propo}

Write
$$
Y'_{P(z)}(v, x)=\tau_{P(z)}(Y_{t}(v, x))
$$
and
$$
Y'_{P(z)}(\omega, x)=\sum_{n\in \Bbb{Z}}L'_{P(z)}(n)x^{-n-2}.
$$
We call the eigenspaces of the operator $L'_{P(z)}(0)$ the {\it weight
subspaces} or {\it homogeneous subspaces} of $(W_{1}\otimes
W_{2})^{*}$, and we have the corresponding notions of {\it weight
vector} (or {\it homogeneous vector}) and {\it weight}.

Let $W$ be a subspace of $(W_{1}\otimes W_{2})^{*}$.  We say that $W$
is {\it compatible for $\tau_{P(z)}$} if every element of $W$
satisfies the following nontrivial and subtle condition (called the
{\it compatibility condition}) on $\lambda
\in (W_{1}\otimes W_{2})^{*}$: The formal Laurent series $Y'_{P(z)}(v,
x_{0})\lambda$ involves only finitely many negative powers of $x_{0}$
and
\begin{eqnarray*}
\lefteqn{\tau_{P(z)}\left(x_{0}^{-1}
\delta\left(\frac{x^{-1}_{1}-z}{x_{0}}
\right)
Y_{t}(v, x_{1})\right)\lambda=}\nno\\
&&=x_{0}^{-1}\delta\left(\frac{x^{-1}_{1}-z}{x_{0}}\right)
Y'_{P(z)}(v, x_{1})\lambda  \;\;\;\;\;
\mbox{\rm for all}\;\;v\in V.\label{comp}
\end{eqnarray*}
Also, we say that $W$ is ($\Bbb{C}$-){\it graded} if it is
$\Bbb{C}$-graded by its weight subspaces, and that $W$ is a $V$-{\it module}
(respectively, {\it generalized module}) if $W$ is graded and is a
module (respectively, generalized module, see \cite{HL1} and
\cite{HL2}) when equipped with this grading and with the action of
$Y'_{P(z)}(\cdot, x)$. The weight subspace of a subspace $W$ with
weight $n\in \Bbb{C}$ will be denoted $W_{(n)}$.

Define
$$
W_{1}\hboxtr_{P(z)}W_{2}=\sum_{W\in {\cal W}_{P(z)}}W =\bigcup_{W\in
{\cal W}_{P(z)}} W\subset
(W_{1}\otimes W_{2})^{*},
$$
where ${\cal W}_{P(z)}$ is the set of all compatible
modules for $\tau_{P(z)}$ in $(W_{1}\otimes W_{2})^{*}$.

We need:

\begin{defi}
{\rm A  vertex operator algebra $V$ is {\it rational} if it
satisfies the following conditions:
\begin{enumerate}
\item There are only finitely many irreducible $V$-modules (up to equivalence).
\item Every $V$-module is completely reducible (and is in particular a
{\it finite} direct sum of irreducible modules).
\item All the fusion rules (the dimensions of spaces of intertwining operators)
for $V$ are finite (for triples of irreducible
modules and hence arbitrary modules).
\end{enumerate}
}
\end{defi}

We have:
\begin{propo}\label{rational}
Let $V$ be a rational vertex operator algebra and $W_{1}$, $W_{2}$
$V$-modules. Then $(W_{1}\hboxtr _{P(z)}W_{2},
\left. Y'_{P(z)}\right|_{V\otimes W_{1}\shboxtr _{P(z)}W_{2}})$
is a module.
\end{propo}

If $W_{1}\hboxtr _{P(z)}W_{2}$ is
a module, we define a
$V$-module $W_{1}\boxtimes_{P(z)} W_{2}$ by
$$
W_{1}\boxtimes_{P(z)} W_{2}=(W_{1}\hboxtr_{P(z)}W_{2})'
$$
 and we write the corresponding action as $Y_{P(z)}$. Applying
Proposition \ref{intwact} to the special module
$W_{3}=W_{1}\boxtimes_{P(z)} W_{2}$ and the identity map $W'_{3}\to
W_{1}\hboxtr_{P(z)} W_{2}$, we obtain  a canonical
$P(z)$-intertwining map of type ${W_{1}\boxtimes_{P(z)} W_{2}}\choose
{W_{1}\;\;W_{2}}$, which we denote
\begin{eqnarray*}
\boxtimes_{P(z)}: W_{1}\otimes W_{2}&\to &
\overline{W_{1}\boxtimes_{P(z)} W_{2}}\nno\\
w_{(1)}\otimes  w_{(2)}&\mapsto& w_{(1)} \boxtimes_{P(z)}w_{(2)}.
\end{eqnarray*}
We have:

\begin{propo}\label{pztensor}
Assume that $W_{1}\hboxtr _{P(z)}W_{2}$ is a module.
Then the $P(z)$-product $(W_{1}\boxtimes_{P(z)} W_{2}, Y_{P(z)};
\boxtimes_{P(z)})$
is a $P(z)$-tensor product of $W_{1}$ and $W_{2}$.
\end{propo}

Combining this proposition with Proposition \ref{rational}, we
obtain:

\begin{propo}
Assume that $V$ is rational. Then $(W_{1}\boxtimes_{P(z)} W_{2}, Y_{P(z)};
\boxtimes_{P(z)})$
is a $P(z)$-tensor product of $W_{1}$ and $W_{2}$.
\end{propo}

Observe that any element of
$W_{1}\hboxtr_{P(z)} W_{2}$ is an element $\lambda$
of $(W_{1}\otimes W_{2})^{*}$ satisfying:

\begin{description}
\item[The compatibility condition] \hfill

{\bf (a)} The  {\it lower
truncation condition}:
For all $v\in V$, the formal Laurent series $Y'_{P(z)}(v, x)\lambda$
involves only finitely many negative
powers of $x$.

{\bf (b)} The formula (\ref{comp}) holds.

\item[The local grading-restriction  condition]\hfill

{\bf (a)} The {\it grading condition}:
$\lambda$ is a (finite) sum of
weight vectors of $(W_{1}\otimes W_{2})^{*}$.

{\bf (b)} Let $W_{\lambda}$ be the smallest subspace of $(W_{1}\otimes
W_{2})^{*}$ containing $\lambda$ and stable under the component
operators $\tau_{P(z)}(v\otimes t^{n})$ of the operators $Y'_{P(z)}(v,
x)$ for $v\in V$, $n\in \Bbb{Z}$. Then the weight spaces
$(W_{\lambda})_{(n)}$, $n\in \Bbb{C}$, of the (graded) space
$W_{\lambda}$ have the properties
\begin{eqnarray*}
&\mbox{\rm dim}\ (W_{\lambda})_{(n)}<\infty \;\;\;\mbox{\rm for}\
n\in \Bbb{C},&\\
&(W_{\lambda})_{(n)}=0 \;\;\;\mbox{\rm for $n$ whose real part is
sufficiently small.}&
\end{eqnarray*}
\end{description}

We have
another construction
of $W_{1}\hboxtr_{P(z)} W_{2}$ using these conditions:

\begin{theo}
The subspace of $(W_{1}\otimes W_{2})^{*}$ consisting of the
elements satisfying the compatibility
condition and the local grading-restriction condition, equipped with
$Y'_{P(z)}$, is a generalized module and is equal to
$W_{1}\hboxtr_{P(z)} W_{2}$.
\end{theo}

The following results follows immediately {from}
 the theorem above, the definition of
$W_{1}\boxtimes_{P(z)} W_{2}$ and Proposition
\ref{rational}:

\begin{corol}\label{pztensor2}
Assume that $W_{1}\hboxtr _{P(z)}W_{2}$ is a module.
Then the contragredient module of the module $W_1
\hboxtr_{P(z)} W_2$,
equipped with the $P(z)$-intertwining map $\boxtimes_{P(z)}$,
 is a $P(z)$-tensor product of $W_{1}$ and
$W_{2}$  equal to the structure $(W_{1}\boxtimes_{P(z)}
W_{2}, Y_{P(z)};
\boxtimes_{P(z)})$ constructed above.
\end{corol}

\begin{corol}
Assume that $V$ is a rational vertex operator algebra.
Then the contragredient module of the module $W_1
\hboxtr_{P(z)} W_2$,
equipped with the $P(z)$-intertwining map $\boxtimes_{P(z)}$,
 is a $P(z)$-tensor product of $W_{1}$ and
$W_{2}$  equal to the structure $(W_{1}\boxtimes_{P(z)}
W_{2}, Y_{P(z)};
\boxtimes_{P(z)})$ constructed above.
\end{corol}

The rationality of a vertex operator algebra is sufficient for the
construction of the tensor product of two modules for this vertex
operator algebra, but is still not enough to
guarantee that the tensor products satisfy a certain associativity
property and that the category of modules is a vertex tensor category defined
in \cite{HL4}. So we need more conditions.

We recall these other conditions. Assume that $V$ is rational and all
irreducible $V$-modules are ${\Bbb R}$-graded. Note that in this case
all $V$-modules are ${\Bbb R}$-graded, that is, weights of elements of
$V$-modules are always real numbers and that all intertwining
operators are formal series of real powers.  Given any $V$-modules
$W_{1}$, $W_{2}$, $W_{3}$, $W_{4}$ and $W_{5}$, let ${\cal Y}_{1}$,
${\cal Y}_{2}$, ${\cal Y}_{3}$ and ${\cal Y}_{4}$ be intertwining
operators of type ${W_{4}}\choose {W_{1}W_{5}}$, ${W_{5}}\choose
{W_{2}W_{3}}$, ${W_{5}}\choose {W_{1}W_{2}}$ and ${W_{4}}\choose
{W_{5}W_{3}}$, respectively. For any $V$-module $W$, let $P_{n}$,
$n\in {\Bbb R}$, be the
projection {}from  $W$ to
its homogeneous $W_{(n)}$subspace of weight $n$.
For any $w_{(1)}\in W_{1}$, $w_{(2)}\in W_{2}$,
$w_{(3)}\in W_{3}$ and $w'_{(4)}\in W'_{4}$, we say that
\begin{equation}\label{1.0}
\langle w'_{(4)}, {\cal Y}_{1}(w_{(1)}, x_{2})
{\cal Y}_{2}(w_{(2)}, x_{2})w_{(3)}\rangle_{W_{4}}
\lbar_{x_{1}= z_{1}, \;x_{2}=z_{2}}
\end{equation}
{\it is convergent} if the series
$$
\sum_{n\in {\Bbb R}}\langle w'_{(4)}, {\cal Y}_{1}(w_{(1)}, x_{2})
P_{n}({\cal Y}_{2}(w_{(2)}, x_{2})w_{(3)})\rangle_{W_{4}}
\lbar_{x_{1}= z_{1}, \;x_{2}=z_{2}}
$$
is absolutely convergent for any choices of $\log z_{1}$ and $\log z_{2}$
in the definitions of $z_{1}^{n}=e^{n\log z_{1}}$ and
$z_{2}^{n}=e^{n\log z_{2}}$, $n\in {\Bbb R}$.
Similarly, we say that
\begin{equation}\label{1.0.5}
\langle w'_{(4)}, {\cal Y}_{4}({\cal Y}_{3}(w_{(1)}, x_{0})
w_{(2)}, x_{2})w_{(3)}\rangle_{W_{4}}
\lbar_{x_{0}= z_{1}-z_{2}, \;x_{2}=z_{2}}
\end{equation}
{\it is convergent} if the series
$$
\sum_{n\in {\Bbb R}}\langle w'_{(4)}, {\cal Y}_{4}(P_{n}({\cal Y}_{3}(w_{(1)},
x_{0})
w_{(2)}), x_{2})w_{(3)}\rangle_{W_{4}}
\lbar_{x_{0}= z_{1}-z_{2}, \;x_{2}=z_{2}}
$$
is absolutely convergent for any choices of $\log (z_{1}-z_{2})$
and $\log z_{2}$
in the definitions of $(z_{1}-z_{2})^{n}=e^{n\log (z_{1}-z_{2})}$ and
$z_{2}^{n}=e^{n\log z_{2}}$, $n\in {\Bbb R}$.
Consider the following conditions for the
product of ${\cal Y}_{1}$ and ${\cal Y}_{2}$ and for the iterate of
${\cal Y}_{3}$ and ${\cal Y}_{4}$, respectively:

\begin{description}

\item[Convergence and extension property for products]
There exists
an integer $N$
(depending only on ${\cal Y}_{1}$ and ${\cal Y}_{2}$), and
for any $w_{(1)}\in W_{1}$,
$w_{(2)}\in W_{2}$, $w_{(3)}\in W_{3}$, $w'_{(4)}\in W'_{4}$, there exist
$j\in {\Bbb N}$, $r_{i}, s_{i}\in {\Bbb R}$, $i=1, \dots, j$, and analytic
functions $f_{i}(z)$ on $|z|<1$, $i=1, \dots, j$,
satisfying
\begin{equation}\label{si}
\wt w_{(1)}+\wt w_{(2)}+s_{i}>N,\;\;\;i=1, \dots, j,
\end{equation}
such that
$$
\langle w'_{(4)}, {\cal Y}_{1}(w_{(1)}, x_{2})
{\cal Y}_{2}(w_{(2)}, x_{2})w_{(3)}\rangle_{W_{4}}
\lbar_{x_{1}= z_{1}, \;x_{2}=z_{2}}
$$
is convergent when $|z_{1}|>|z_{2}|>0$ and can be analytically extended to
the multi-valued analytic function
\begin{equation}\label{phyper}
\sum_{i=1}^{j}z_{2}^{r_{i}}(z_{1}-z_{2})^{s_{i}}
f_{i}\left(\frac{z_{1}-z_{2}}{z_{2}}\right)
\end{equation}
when $|z_{2}|>|z_{1}-z_{2}|>0$.

\item[Convergence and extension property for iterates] There exists an integer
$\tilde{N}$
(depending only on ${\cal Y}_{3}$ and ${\cal Y}_{4}$), and
for any $w_{(1)}\in W_{1}$,
$w_{(2)}\in W_{2}$, $w_{(3)}\in W_{3}$, $w'_{(4)}\in W'_{4}$, there exist
$k\in {\Bbb N}$, $\tilde{r}_{i}, \tilde{s}_{i}\in {\Bbb R}$, $i=1, \dots, k$,
and analytic
functions $\tilde{f}_{i}(z)$ on $|z|<1$, $i=1, \dots, k$,
satisfying
$$
\wt w_{(2)}+\wt w_{(3)}+\tilde{s}_{i}>\tilde{N},\;\;\;i=1, \dots, k,
$$
 such that
$$
\langle w'_{(4)},
{\cal Y}_{4}({\cal Y}_{3}(w_{(1)}, x_{0})w_{(2)}, x_{2})w_{(3)}\rangle_{W_{4}}
\lbar_{x_{0}=z_{1}-z_{2},\;x_{2}=z_{2}}
$$
is convergent when $|z_{2}|>|z_{1}-z_{2}|>0$ and can be analytically extended
to the multi-valued analytic function
$$
\sum_{i=1}^{k}z_{1}^{\tilde{r}_{i}}z_{2}^{\tilde{s}_{i}}
\tilde{f}_{i}\left(\frac{z_{2}}{z_{1}}\right)
$$
when $|z_{1}|>|z_{2}|>0$.
\end{description}

If for any $V$-modules $W_{1}$, $W_{2}$, $W_{3}$, $W_{4}$ and $W_{5}$ and
 any intertwining operators ${\cal Y}_{1}$ and ${\cal Y}_{2}$
of the types as above, the convergence and extension property for products
holds,
we say that
{\it  the products of the
 intertwining operators for $V$ have the convergence and extension property}.
Similarly we can define the meaning of the phrase
{\it the iterates of the intertwining
operators for $V$ have the  convergence and extension property}.

We also need a condition on certain generalized modules.
If a generalized $V$-module $W=\coprod_{n\in \Bbb{C}}W_{(n)}$
satisfying the condition that $W_{(n)}=0$ for $n$ whose real part is
sufficiently small, we say that $W$ is {\it lower-truncated}.
Another condition that we need to establish the associativity
is that every finitely-generated lower-truncated
generalized $V$-module is a $V$-module.

Assume that for any $V$-modules $W_{i}$, $i=1, \dots, 5$,
any intertwining operators ${\cal Y}_{1}$, ${\cal Y}_{2}$
of the type above, (\ref{1.0}) is convergent
for all $w_{(1)}\in W_{1}$, $w_{(2)}\in W_{2}$,
$w_{(3)}\in W_{3}$ and $w'_{(4)}\in W'_{4}$, and
for any $V$-modules $W_{i}$, $i=1, \dots, 5$,
any intertwining operators ${\cal Y}_{3}$, ${\cal Y}_{4}$
of the type above, (\ref{1.0.5}) is convergent
for all $w_{(1)}\in W_{1}$, $w_{(2)}\in W_{2}$,
$w_{(3)}\in W_{3}$ and $w'_{(4)}\in W'_{4}$.
Let $W_{1}$, $W_{2}$ and $W_{3}$ be
three $V$-modules, $w_{(1)}\in W_{1}$, $w_{(2)}\in W_{2}$ and
$w_{(3)}\in W_{3}$ and $z_{1}, z_{2}\in \Bbb{C}$ satisfying
$|z_{1}|>|z_{2}|>|z_{1}-z_{2}|>0$. {}From \cite{HL2},
we know that any $P(z)$-intertwining maps (for
$z=z_{1}, z_{2}, z_{1}-z_{2}$) can be obtained {}from
certain intertwining
operators by substituting complex powers of $e^{\log z}$ for the
complex powers of the formal variable $x$. Thus
$w_{(1)}\boxtimes_{P(z_{1})}(w_{(2)}\boxtimes_{P(z_{2})}w_{(3)})$ (or
$(w_{(1)}\boxtimes_{P(z_{1}-z_{2})}w_{(2)})\boxtimes_{P(z_{2})}w_{(3)}$)
is a product (or an iterate) of two intertwining operators evaluated
at $w_{(1)}\otimes w_{(2)}\otimes w_{(3)}$ and with the complex powers
of the formal variables replaced by the complex powers of $e^{\log
z_{1}}$ and of $e^{\log z_{2}}$ (or by the complex powers of $e^{\log
(z_{1}-z_{2})}$ and of $e^{\log z_{2}}$). By assumption,
$w_{(1)}\boxtimes_{P(z_{1})}(w_{(2)}\boxtimes_{P(z_{2})}w_{(3)})$ (or
$(w_{(1)}\boxtimes_{P(z_{1}-z_{2})}w_{(2)})\boxtimes_{P(z_{2})}w_{(3)}$)
 is a well-defined element of
$\overline{W_{1}\boxtimes_{P(z_{1})}(W_{2}\boxtimes_{P(z_{2})}W_{3})}$
(or of $\overline{(W_{1}\boxtimes_{P(z_{1}-z_{2})}W_{2})
\boxtimes_{P(z_{2})}W_{3}}\;$).
The following associativity of the $P(\cdot)$-tensor products
is proved in \cite{H1}:

\begin{theo}\label{assoc}
Assume that $V$ is a rational vertex operator algebra and all irreducible
$V$-modules are ${\Bbb R}$-graded.  Also assume
that $V$ satisfies the following conditions:
\begin{enumerate}
\item Every finitely-generated lower-truncated generalized $V$-module
is a $V$-module.

\item The products or the iterates of
the intertwining operators for $V$ have the  convergence and
extension property.
\end{enumerate}
Then for any $V$-module $W_{1}$, $W_{2}$ and $W_{3}$ and any complex
numbers $z_{1}$ and $z_{2}$ satisfying
$|z_{1}|>|z_{2}|>|z_{1}-z_{2}|>0$, there exists a unique isomorphism
${\cal A}_{P(z_{1}), P(z_{2})}^{P(z_{1}-z_{2}), P(z_{2})}$ {}from
$W_1\boxtimes_{P(z_1)}(W_2 \boxtimes_{P(z_2)}W_3)$ to $(W_1
\boxtimes_{P(z_1-z_2)}W_2)\boxtimes_{P(z_2)}W_3$ such that
for any $w_{(1)}\in W_{1}$, $w_{(2)}\in W_{2}$ and $w_{(3)}\in W_{3}$,
\begin{eqnarray*}
\lefteqn{\overline{{\cal A}}_{P(z_{1}), P(z_{2})}^{P(z_{1}-z_{2}),
P(z_{2})}(w_{(1)}\boxtimes_{P(z_{1})}(w_{(2)}
\boxtimes_{P(z_{2})}w_{(3)}))}\nno\\
&&=(w_{(1)}\boxtimes_{P(z_{1}-z_{2})}w_{(2)})
\boxtimes_{P(z_{2})}w_{(3)},
\end{eqnarray*}
where $$\overline{{\cal A}}_{P(z_{1}), P(z_{2})}^{P(z_{1}-z_{2}),
P(z_{2})}: \overline{W_1\boxtimes_{P(z_1)}(W_2 \boxtimes_{P(z_2)}W_3)}
\to \overline{(W_1 \boxtimes_{P(z_1-z_2)}W_2)\boxtimes_{P(z_2)}W_3}$$
is the unique extension of ${\cal A}_{P(z_{1}),
P(z_{2})}^{P(z_{1}-z_{2}), P(z_{2})}$.
\end{theo}

In \cite{H1}, it is also proved that the above associativity of the
$P(\cdot)$-tensor product is equivalent to the associativity of the
intertwining operators. So we also have the following:

\begin{theo}\label{intwassoc}
Let $V$ be a vertex operator algebra satisfying the conditions
in Theorem \ref{assoc}.
Then the intertwining operators for
$V$ have the following two
associativity properties:
\begin{enumerate}

\item For any modules
$W_{1}$, $W_{2}$, $W_{3}$, $W_{4}$ and $W_{5}$ and
any
intertwining operators ${\cal Y}_{1}$ and ${\cal Y}_{2}$ of
 type ${W_{4}}\choose {W_{1}W_{5}}$ and ${W_{5}}\choose {W_{2}W_{3}}$,
respectively,
there exist a module $W_{6}$ and intertwining operators ${\cal Y}_{3}$
and ${\cal Y}_{4}$ of type ${W_{6}}\choose {W_{1}W_{2}}$ and
${W_{4}}\choose {W_{6}W_{3}}$, respectively, such that
for any $z_{1}, z_{2}\in {\Bbb C}$
satisfying $|z_{1}|>|z_{2}|>|z_{1}-z_{2}|>0$,
\begin{eqnarray}\label{1-1}
\lefteqn{\langle w'_{(4)}, {\cal Y}_{1}(w_{(1)},
x_{1}){\cal Y}_{2}(w_{(2)}, x_{2})w_{(3)}\rangle_{W_{4}}\lbar_{x^{n}_{1}
=e^{n\log z_{1}},\;
x^{n}_{2}=e^{n\log z_{2}}, \; n\in {\Bbb C}}}
\nno\\
&&=\langle w'_{(4)}, {\cal Y}_{4}({\cal Y}_{3}(w_{(1)},
x_{0})w_{(2)}, x_{2})w_{(3)}\rangle_{W_{4}}\lbar_{x^{n}_{0}=e^{n \log
(z_{1}-z_{2})},\;
x^{n}_{2}=e^{n\log z_{2}}, \; n\in {\Bbb C}}.\nno\\
&&
\end{eqnarray}
 holds for any
$w'_{(4)}\in W'_{4}$, $w_{(1)}\in W_{1}$,
$w_{(2)}\in W_{2}$ and $w_{(3)}\in W_{3}$.

\item For any modules
$W_{1}$, $W_{2}$, $W_{3}$, $W_{4}$ and $W_{6}$ and
any intertwining operators ${\cal Y}_{3}$
and ${\cal Y}_{4}$ of type ${W_{6}}\choose {W_{1}W_{2}}$ and
${W_{4}}\choose {W_{6}W_{3}}$, respectively, there exist a module $W_{5}$
and intertwining operators ${\cal Y}_{1}$ and ${\cal Y}_{2}$ of
 type ${W_{4}}\choose {W_{1}W_{5}}$ and ${W_{5}}\choose {W_{2}W_{3}}$,
respectively, such that
for any $z_{1}, z_{2}\in {\Bbb C}$
satisfying $|z_{1}|>|z_{2}|>|z_{1}-z_{2}|>0$, (\ref{1-1}) holds for any
$w'_{(4)}\in W'_{4}$, $w_{(1)}\in W_{1}$,
$w_{(2)}\in W_{2}$ and $w_{(3)}\in W_{3}$.

\end{enumerate}
\end{theo}

The result below is a consequence of the main result announced in
\cite{HL4} and Theorem \ref{intwassoc}. It
gives the structure of a vertex tensor category to
the category of modules for a vertex operator algebra satisfying the conditions
discussed above and an additional convergence condition.
See \cite{HL4} for the definitions of
$\tilde{K}^{c}$, $\psi$ and vertex tensor
category.

\begin{theo}\label{vtc}
Let $V$ be a vertex operator algebra of central charge $c$
satisfying the conditions
in Theorem \ref{assoc}. We assume in addition that $V$ satisfies
the following condition:
\begin{enumerate}

\setcounter{enumi}{2}

\item  For any $V$-modules $W_{j}$, $j=1, \dots, 2m+1$, any
intertwining operators ${\cal  Y}_{i}$, $i=1, \dots, m$,
of types ${W_{2i-1}}\choose {W_{2i}W_{2i+1}}$, respectively,
and any $w'_{(1)}\in W'_{1}$, $w_{(2i)}\in W_{2i}$, $i=1, \dots, m$,
and $w_{(2m+1)}\in W_{2m+1}$,
$$\langle w'_{1}, {\cal  Y}_{1}(w_{(2)}, x_{1})\cdots
{\cal  Y}_{m}(w_{(2m)}, x_{m})w_{(2m+1)}\rangle\lbar_{x_{i}^{n}
=e^{n\log z_{i}}, \; 1\le i\le m, \; n\in {\Bbb C}}$$
is absolutely convergent
 for any $z_{1}, \dots, z_{n}\in {\Bbb C}$
satisfying $|z_{1}|>\cdots >|z_{n}|>0$.
\end{enumerate}
Then the category of $V$-modules has a natural structure of a vertex
tensor category of central charge $c$ such that for each $z\in {\Bbb
C}^{\times}$, the tensor product bifunctor $\boxtimes_{\psi(P(z))}$
associated to $\psi(P(z))\in
\tilde{K}^{c}(2)$ is equal to $\boxtimes_{P(z)}$.
\end{theo}

The following result is also announced in \cite{HL4} (see \cite{JS} for the
definition of braided tensor (monoidal) category):

\begin{theo}\label{btc}
The underlying category of a vertex tensor category has a natural
structure of a braided tensor category.
\end{theo}

\renewcommand{\theequation}{\thesection.\arabic{equation}}
\renewcommand{\therema}{\thesection.\arabic{rema}}
\setcounter{equation}{0}
\setcounter{rema}{0}

\section{Minimal Virasoro vertex operator algebras and modules}

In this section, we summarize the results on minimal Virasoro vertex
operator algebras obtained in \cite{FZ}, \cite{DMZ} and \cite{W}.

Let
$$
{\frak L}=\oplus_{n\in {\Bbb Z}}{\Bbb C}L_{n}\oplus {\Bbb C}d
$$
with the commutation relations
\begin{eqnarray*}
{[L_{m}, L_{n}]}&=&(m-n)L_{m+n}+\frac{m^{3}-m}{12}\delta_{m+n, 0}d,\;\;\;
m, n\in {\Bbb Z}\\
{[L_{m}, d]}&=&0, \;\;\;m\in {\Bbb Z},
\end{eqnarray*}
 be the Virasoro algebra. Consider the two subalgebras
\begin{eqnarray*}
{\frak L}_{+}&=&\oplus_{n\in {\Bbb Z}_{+}}{\Bbb C}L_{n},\\
{\frak L}_{-}&=&\oplus_{n\in -{\Bbb Z}_{+}}{\Bbb C}L_{n}.
\end{eqnarray*}
Let $U(\cdot)$ be the functor {}from the category of Lie algebras to
the category of associative algebras obtained by taking the universal
enveloping algebras of  Lie algebras.
For any representation of ${\frak L}$ and any
$m\in {\Bbb Z}$, we shall use
$L(m)$ to denote the representation image of $L_{m}$.
For any $h, c\in {\Bbb C}$, the Verma module $M(c, h)$ for ${\frak L}$
is a free $U({\frak L}_{-})$-module generated by ${\bf 1}_{c, h}$ such that
\begin{eqnarray*}
{\frak L}_{+}{\bf 1}_{c, h}&=&0,\\
L(0){\bf 1}_{c, h}&=&h{\bf 1}_{c, h},\\
d{\bf 1}_{c, h}&=&c{\bf 1}_{c, h}.
\end{eqnarray*}
There exists a unique maximal proper submodule $J(c, h)$ of $M(c, h)$. It
is easy to see that if $c\ne 0$, then
both ${\bf 1}_{c, 0}$ and $L(-2){\bf 1}_{c, 0}$
are not in $J(c, 0)$.
Let
$$
L(c, h)=M(c, h)/J(c, h).
$$
Then $L(c, 0)$ has the structure of a vertex operator algebra with
vacuum ${\bf 1}_{c, 0}$ and the Virasoro element $L(-2){\bf 1}_{c, 0}$
\cite{FZ}.  The following result is proved in \cite{W} using the
results of Feigin and Fuchs on representations of the Virasoro algebra
\cite{FF1} \cite{FF2}:

\begin{theo}
The vertex operator algebra $L(c, 0)$ is rational if and only if either $c=0$
or
there is a pair $p, q$ of relatively prime positive integers larger
than $1$ such that
$$
c=c_{p, q}=1-6\frac{(p-q)^{2}}{pq}.
$$
A set of representatives of the equivalence classes of irreducible modules
for $L(c_{p, q}, 0)$ is
$$
\{L(c_{p, q}, h_{p, q}^{m, n})\}_{0<m<p, \; 0<n<q, \; m, n\in {\Bbb Z}}
$$
where for any $m, n\in {\Bbb Z}$ satisfying $0<m<p, 0<n<q$,
$$
h_{p, q}^{m, n}=\frac{(np-mq)^{2}-(p-q)^{2}}{4pq}.
$$
For any $m, m', m'', n, n', n''\in {\Bbb Z}$ satisfying
$0<m, m', m''<p, 0<n, n', n''<q$, the fusion rule ${\cal N}_{L(c_{p, q}, \;
h_{p, q}^{m', n'})L(c_{p, q}, \;
h_{p, q}^{m'', n''})}^{L(c_{p, q}, \;
h_{p, q}^{m, n})}$ is $1$ if $m+m'+m''<2p$, $n+n'+n''<2q$, $m<m'+m''$,
$m'<m''+m$, $m''<m+m'$, $n<n'+n''$, $n'<n''+n$, $n''<n+n'$ and the
sums $m+m'+m''$, $n+n'+n''$ are odd, and is $0$ otherwise.
\end{theo}

The rationality of $L(c_{p, p+1}, 0)$ for an integer $p>1$
was proved first in \cite{DMZ}
and the fusion rules in the case $p=3$, $q=4$ were also calculated
there.

For any pair $p, q$ of relatively prime positive integers larger
than $1$, we call the vertex operator algebra $L(c_{p, q}, 0)$ a
{\it minimal Virasoro vertex operator algebra}.

Let $n$ be a positive integer, $(p_{i}, q_{i})$, $i=1, \dots, n$, $n$
pairs of relatively prime positive integers larger than $1$,
$V=L(c_{p_{1}, q_{1}}, 0)\otimes \cdots \otimes L(c_{p_{n}, q_{n}},
0)$. {}From the results proved in \cite{FHL} and \cite{DMZ}, $V$ is a
rational vertex operator algebra, a set of representatives of equivalence
classes of irreducible modules for $V$ can be given explicitly and the
fusion rules for $V$ can be calculated easily.

\renewcommand{\theequation}{\thesection.\arabic{equation}}
\renewcommand{\therema}{\thesection.\arabic{rema}}
\setcounter{equation}{0}
\setcounter{rema}{0}

\section{The results of the present paper}

First we have the following commutativity for intertwining operators:

\begin{theo}\label{commu}
Let $V$ be a vertex operator algebra satisfying the conditions
in Theorem \ref{assoc}. Then for any modules
$W_{1}$, $W_{2}$, $W_{3}$, $W_{4}$ and $W_{5}$ and
any
intertwining operators ${\cal Y}_{1}$ and ${\cal Y}_{2}$ of
 type ${W_{4}}\choose {W_{1}W_{5}}$ and ${W_{5}}\choose {W_{2}W_{3}}$,
respectively,
there exist a module $W_{6}$ and intertwining operators ${\cal Y}_{3}$
and ${\cal Y}_{4}$ of type ${W_{4}}\choose {W_{2}W_{6}}$ and
${W_{6}}\choose {W_{1}W_{3}}$, respectively, such that for any
$w'_{(4)}\in W'_{4}$, $w_{(1)}\in W_{1}$,
$w_{(2)}\in W_{2}$ and $w_{(3)}\in W_{3}$,
the multi-valued analytic function
$$\langle w'_{(4)}, {\cal Y}_{1}(w_{(1)},
x_{1}){\cal Y}_{2}(w_{(2)}, x_{2})w_{(3)}\rangle_{W_{4}}\lbar_{x_{1}
=z_{1},\;
x_{2}=z_{2}}$$
of $z_{1}$ and $z_{2}$ in the region $|z_{1}|>|z_{2}|>0$ and
the multi-valued analytic function
$$\langle w'_{(4)}, {\cal Y}_{3}(w_{(2)},
x_{2}){\cal Y}_{4}(w_{(1)}, x_{1})w_{(3)}\rangle_{W_{4}}\lbar_{x_{1}
=z_{1},\;
x_{2}=z_{2}}$$
of $z_{1}$ and $z_{2}$ in the region $|z_{2}|>|z_{1}|>0$, are analytic
extensions of each other.
\end{theo}

The next result relaxes the conditions to use the tensor product theory:

\begin{theo}\label{subvtc}
Let $V$ be a vertex operator algebra containing a rational
vertex operator subalgebra
$V_{0}$.
Then for any $V$-modules $W_{1}$ and $W_{2}$,
$W_{1}\hboxtr_{P(z)}W_{2}$ is a $V$-module and
the conclusions of Proposition \ref{pztensor} and
Corollary \ref{pztensor2} are true.
Suppose in addition that we assume
that all irreducible
$V_{0}$-modules are ${\Bbb R}$-graded and $V$ and
$V_{0}$ satisfies the following conditions:
\begin{enumerate}
\item Every finitely-generated lower-truncated generalized $V$-module
is a $V$-module.

\item The products or the iterates of
the intertwining operators for $V_{0}$ have the  convergence and
extension property.
\end{enumerate}
Then
the conclusions of  Theorem \ref{assoc}, Theorem \ref{intwassoc}
and Theorem \ref{commu} are true. Suppose that
 we assume further that $V_{0}$ satisfies
the following condition:

\begin{enumerate}

\setcounter{enumi}{2}

\item For any $V_{0}$-modules $W_{j}$, $j=1, \dots, 2m+1$, any
intertwining operators ${\cal  Y}_{i}$, $i=1, \dots, m$,
of types ${W_{2i-1}}\choose {W_{2i}W_{2i+1}}$, respectively,
and any $w'_{(1)}\in W'_{1}$, $w_{(2i)}\in W_{2i}$, $i=1, \dots, m$,
and $w_{(2m+1)}\in W_{2m+1}$,
$$\langle w'_{1}, {\cal  Y}_{1}(w_{(2)}, x_{1})\cdots
{\cal  Y}_{m}(w_{(2m)}, x_{m})w_{(2m+1)}\rangle\lbar_{x_{i}^{n}
=e^{n\log z_{i}}, \; 1\le i\le m, \; n\in {\Bbb C}}$$
is absolutely convergent
 for any $z_{1}, \dots, z_{n}\in {\Bbb C}$
satisfying $|z_{1}|>\cdots >|z_{n}|>0$.
\end{enumerate}
Then the conclusion of Theorem \ref{vtc} is true.
\end{theo}

By Theorem \ref{subvtc} and Theorem \ref{btc}, we obtain:

\begin{corol}\label{subbtc}
Let $V$ be a vertex operator algebra satisfying all conditions in
Theorem \ref{subvtc}. Then the category of $V$-modules has a natural
structure of a braided tensor category.
\end{corol}

Because of the above theorem,
It is natural to ask whether a vertex operator algebra containing a
rational vertex operator subalgebra is also rational. The answer is no.
Here is a simple example:

\begin{exam}\label{irra}
{\rm Let $(V_{0}, Y_{V_{0}}, {\bf 1}, \omega)$ be a rational vertex operator
algebra.
Let $V=V_{0}\oplus V_{0}$. Then $V$ is ${\Bbb Z}$-graded and satisfies
the two grading axioms for vertex operator algebras. We denote
an element of $V$ by $(u, v)$ where $u, v\in V_{0}$. We define a vertex
operator map $Y_{V}: V\otimes V\to V[[x, x^{-1}]]$ by
$$
Y_{V}((u_{1}, v_{1}), x)(u_{2}, v_{2})=
(Y_{V_{0}}(u_{1}, x)u_{2}, Y_{V_{0}}(u_{1}, x)v_{2}+
Y_{V_{0}}(v_{1}, x)u_{2})
$$
for any $(u_{1}, v_{1}), (u_{2}, v_{2})\in V$.
Then using the vertex operator algebra structure on $V_{0}$,
we see that
 $(V, Y_{V}, ({\bf 1}, 0), (\omega, 0))$
is a vertex operator algebra containing
a vertex operator subalgebra isomorphic to $V_{0}$.
Let $W_{0}$ be the subspace of $V$ consists of all
elements of the form $(0, v)$, $v\in V_{0}$.
We define a vertex operator map $Y_{W_{0}}: V\otimes W_{0}\to
W_{0}[[x, x^{-1}]]$ to be the restriction of $Y_{V}$ to $V\otimes W_{0}$.
Then it is clear that $(W_{0}, Y_{W_{0}})$ is a $V$-module.
 The adjoint $V$-module $(V, Y_{V})$ containing
the proper nontrivial $V$-submodule $(W_{0}, Y_{W_{0}})$.
But $(V, Y_{V})$
is not completely reducible since for any element $(u, v)\in V$
such that $u\ne 0$,
$Y_{V}((0, {\bf 1}), x)(u, v)=(0, u)\ne 0$ and is in $W_{0}$.}
\end{exam}

The next result shows that the minimal Virasoro vertex operator algebras
and their tensor product algebras
satisfy, besides the rationality and the rationality of the gradings of
irreducible modules,
also the other conditions to use the tensor product theory:

\begin{theo}\label{cpqm}
For any positive integer $m$ and any
$n$ pairs $(p_{i}, q_{i})$ of relatively prime positive
integers larger than $1$, $i=1, \dots, m$, we have:

\begin{enumerate}

\item Every finitely-generated lower-truncated generalized
$L(c_{p_{1}, q_{1}}, 0)\otimes \cdots \otimes
L(c_{p_{m}, q_{m}}, 0)$-module is a module.

\item The products of the intertwining operators for
$$L(c_{p_{1}, q_{1}}, 0)\otimes \cdots \otimes
L(c_{p_{m}, q_{m}}, 0)$$ have the convergence and extension property.

\item For any modules $W_{j}$, $j=1, \dots, 2m+1$,
for
$$L(c_{p_{1}, q_{1}}, 0)\otimes \cdots \otimes
L(c_{p_{m}, q_{m}}, 0),$$ any
intertwining operators ${\cal  Y}_{i}$, $i=1, \dots, m$,
of types ${W_{2i-1}}\choose {W_{2i}W_{2i+1}}$, respectively,
and any $w'_{(1)}\in W'_{1}$, $w_{(2i)}\in W_{2i}$, $i=1, \dots, m$,
and $w_{(2m+1)}\in W_{2m+1}$,
$$\langle w'_{1}, {\cal  Y}_{1}(w_{(2)}, x_{1})\cdots
{\cal  Y}_{m}(w_{(2m)}, x_{m})w_{(2m+1)}\rangle\lbar_{x_{i}^{n}
=e^{n\log z_{i}}, \; 1\le i\le m, \; n\in {\Bbb C}}$$
is absolutely convergent
 for any $z_{1}, \dots, z_{n}\in {\Bbb C}$
satisfying $|z_{1}|>\cdots >|z_{n}|>0$.

\end{enumerate}
\end{theo}

We now define precisely the
class of vertex operator algebras that we  study in this
paper.

\begin{defi}
{\rm Let $m$ be a nonnegative integer,
$(p_{i}, q_{i})$, $i=1, \dots, m$, $m$ pairs of relatively prime
positive integers larger than $1$.  A vertex operator algebra $V$ is
said to be {\it in the class ${\cal C}_{p_{1}, q_{1};\dots; p_{m}, q_{m}}$}
if $V$ has a vertex operator subalgebra isomorphic to $L(c_{p_{1}, q_{1}},
0)\otimes \cdots \otimes L(c_{p_{n}, q_{n}}, 0)$.}
\end{defi}

We have:

\begin{propo}\label{cpqm1}
Let $V$ be a vertex operator algebra in the class
${\cal C}_{p_{1}, q_{1};\dots; p_{m}, q_{m}}$. Then
every finitely-generated lower-truncated generalized $V$-module
is a $V$-module.
\end{propo}

Combining Theorem \ref{cpqm} and Proposition \ref{cpqm1} with Theorem
\ref{subvtc},
we obtain the main result of the present paper:

\begin{theo}\label{main}
Let $V$ be a vertex operator algebra in the class
${\cal C}_{p_{1}, q_{1};\dots; p_{m}, q_{m}}$.
Then  for any $V$-modules $W_{1}$ and $W_{2}$,
$W_{1}\hboxtr_{P(z)}W_{2}$ is a module and
the conclusions of Proposition \ref{pztensor},
Corollary \ref{pztensor2}, Theorem \ref{assoc}, Theorem \ref{intwassoc},
Theorem \ref{commu}
and Theorem \ref{vtc} are true.
\end{theo}

Combining this main result with Corollary \ref{subbtc}, we obtain:

\begin{corol}\label{cpqmbtc}
Let $V$ be a vertex operator algebra in the class
${\cal C}_{p_{1}, q_{1};\dots; p_{m}, q_{m}}$.
Then the category of $V$-modules has a natural
structure of a braided tensor category. In particular,
the category of $L(c_{p_{1}, q_{1}}, 0)\otimes \cdots \otimes
L(c_{p_{m}, q_{m}}, 0)$-modules has a natural
structure of  a braided tensor category.
\end{corol}

Let $(p, q)$ be a pair of relatively prime
positive integers larger than $1$. We define the {\it category
generated by minimal modules of central charge $c_{p, q}$
for the Virasoro algebra} to be the subcategory of the category of
modules for the
Virasoro algebra such that any object in this subcategory is isomorphic
to a finite direct sum of $L(c_{p, q}, h_{p, q}^{m, n})$, $m, n\in {\Bbb Z}$,
$0<m<p, 0<n<q$. Then the special case $V=L(c_{p, q}, 0)$ in Corollary
\ref{cpqmbtc} can be reformulated as the following:

\begin{theo}\label{virbtc}
Let $(p, q)$ be a pair of relatively prime
positive integers larger than $1$. Then the category
generated by minimal modules of central charge $c_{p, q}$
for the Virasoro algebra has a natural
structure of  a braided tensor category.
\end{theo}

\renewcommand{\theequation}{\thesection.\arabic{equation}}
\renewcommand{\therema}{\thesection.\arabic{rema}}
\setcounter{equation}{0}
\setcounter{rema}{0}

\section{Proof of Theorem 3.1}

By Theorem \ref{intwassoc}, there exist a $V$-module $W_{7}$ and
intertwining operators ${\cal Y}_{5}$ and ${\cal Y}_{6}$ of type
${W_{7}}\choose {W_{1}W_{2}}$ and ${W_{4}}\choose {W_{7}W_{3}}$ such that
for any $w'_{(4)}\in W'_{4}$, $w_{(1)}\in W_{1}$,
$w_{(2)}\in W_{2}$ and $w_{(3)}\in W_{3}$ and
for any $z_{1}, z_{2}\in {\Bbb C}$
satisfying $|z_{1}|>|z_{2}|>|z_{1}-z_{2}|>0$,
\begin{eqnarray}\label{4-1}
\lefteqn{\langle w'_{(4)}, {\cal Y}_{1}(w_{(1)},
x_{1}){\cal Y}_{2}(w_{(2)}, x_{2})w_{(3)}\rangle_{W_{4}}\lbar_{x^{n}_{1}
=e^{n\log z_{1}},\;
x^{n}_{2}=e^{n\log z_{2}}, \; n\in {\Bbb C}}}
\nno\\
&&=\langle w'_{(4)}, {\cal Y}_{6}({\cal Y}_{5}(w_{(1)},
x_{0})w_{(2)}, x_{2})w_{(3)}\rangle_{W_{4}}\lbar_{x^{n}_{0}=e^{n \log
(z_{1}-z_{2})},\;
x^{n}_{2}=e^{n\log z_{2}}, \; n\in {\Bbb C}}.\nno\\
&&
\end{eqnarray}
Substituting
\begin{eqnarray*}
{\cal Y}_{5}(w_{(1)},
x_{0})w_{(2)}&=&\Omega_{0}(\Omega_{1}({\cal Y}_{5}))(w_{(1)},
x_{0})w_{(2)}\nno\\
&=&e^{x_{0}L(-1)}\Omega_{1}({\cal Y}_{5})(w_{(2)}, e^{\pi i}x_{0})w_{(1)}.
\end{eqnarray*}
into (\ref{4-1}), we obtain
\begin{eqnarray}\label{4-2}
\lefteqn{\langle w'_{(4)}, {\cal Y}_{1}(w_{(1)},
x_{1}){\cal Y}_{2}(w_{(2)}, x_{2})w_{(3)}\rangle_{W_{4}}\lbar_{x^{n}_{1}
=e^{n\log z_{1}},\;
x^{n}_{2}=e^{n\log z_{2}}, \; n\in {\Bbb C}}}
\nno\\
&&=\langle w'_{(4)}, {\cal Y}_{6}(e^{x_{0}L(-1)}
\Omega_{1}({\cal Y}_{5})(w_{(2)}, e^{\pi i}x_{0})w_{(1)},
x_{2})\cdot\nno\\
&&\hspace{9em}\cdot w_{(3)}\rangle_{W_{4}}\lbar_{x^{n}_{0}=e^{n \log
(z_{1}-z_{2})},\;
x^{n}_{2}=e^{n\log z_{2}}, \; n\in {\Bbb C}}\nno\\
&&=\langle w'_{(4)}, {\cal Y}_{6}(
\Omega_{1}({\cal Y}_{5})(w_{(2)}, e^{\pi i}x_{0})w_{(1)},
x_{2}+x_{0})\cdot\nno\\
&&\hspace{9em}\cdot w_{(3)}\rangle_{W_{4}}\lbar_{x^{n}_{0}=e^{n \log
(z_{1}-z_{2})},\;
x^{n}_{2}=e^{n\log z_{2}}, \; n\in {\Bbb C}}
\end{eqnarray}
when $|z_{1}|>|z_{2}|>|z_{1}-z_{2}|>0$. The right-hand side of (\ref{4-2})
is the analytic extension of
$$
\langle w'_{(4)}, {\cal Y}_{6}(
\Omega_{1}({\cal Y}_{5})(w_{(2)}, x_{0})w_{(1)},
x_{1})w_{(3)}\rangle_{W_{4}}\lbar_{x^{n}_{0}=e^{n \log
(z_{2}-z_{1})},\;
x^{n}_{1}=e^{n\log z_{1}}, \; n\in {\Bbb C}}
$$
defined in the region $|z_{1}|>|z_{1}-z_{1}|>0$. By Theorem \ref{intwassoc},
there exist a $V$-module $W_{6}$ and intertwining operators ${\cal Y}_{3}$
and ${\cal Y}_{4}$ of type ${W_{4}}\choose {W_{2}W_{6}}$ and
${W_{6}}\choose {W_{1}W_{3}}$, respectively, such that
for any $z_{1}, z_{2}\in {\Bbb C}$
satisfying $|z_{2}|>|z_{1}|>|z_{1}-z_{2}|>0$,
\begin{eqnarray}\label{4-3}
\lefteqn{\langle w'_{(4)}, {\cal Y}_{6}(
\Omega_{1}({\cal Y}_{5})(w_{(2)}, x_{0})w_{(1)},
x_{1})w_{(3)}\rangle_{W_{4}}\lbar_{x^{n}_{0}=e^{n \log
(z_{2}-z_{1})},\;
x^{n}_{1}=e^{n\log z_{1}}, \; n\in {\Bbb C}}}
\nno\\
&&=\langle w'_{(4)}, {\cal Y}_{3}(w_{(2)},
x_{2}){\cal Y}_{4}(w_{(1)}, x_{1})w_{(3)}\rangle_{W_{4}}\lbar_{x^{n}_{1}
=e^{n\log z_{1}},\;
x^{n}_{2}=e^{n\log z_{2}}, \; n\in {\Bbb C}}.
\end{eqnarray}
Thus the left-hand side of (\ref{4-1}) and the right-hand side of
(\ref{4-3}) are analytic extensions of each other, proving the theorem.

\renewcommand{\theequation}{\thesection.\arabic{equation}}
\renewcommand{\therema}{\thesection.\arabic{rema}}
\setcounter{equation}{0}
\setcounter{rema}{0}

\section{Proof of Theorem 3.2}

Let $V$ be a vertex operator algebra containing a rational
vertex operator subalgebra
 $V_{0}$.
We shall denote the $P(z)$-tensor product
operations for $V$ and for $V_{0}$ by $\boxtimes_{P(z)}$ and
$\boxtimes_{P(z)}^{V_{0}}$, respectively. Similarly, we have the notations
$\hboxtr_{P(z)}$ and $\hboxtr_{P(z)}^{V_{0}}$.
Let $W_{1}$ and $W_{2}$ be two $V$-module. Then $W_{1}$ and $W_{2}$ are also
$V_{0}$-modules and by definition we have $ W_{1}\hboxtr_{P(z)}W_{2}
\subset W_{1}\hboxtr^{V_{0}}_{P(z)}W_{2}$.
Since $V_{0}$ is rational, $W_{1}\hboxtr^{V_{0}}_{P(z)}W_{2}$ is
a $V_{0}$-module. In particular, its homogeneous subspaces are
finite-dimensional and when the real part of the weight of a homogeneous
subspace is sufficiently small, the homogeneous subspace
is $0$. Thus $W_{1}\hboxtr_{P(z)}W_{2}$
also satisfies these two conditions. This shows that
$W_{1}\hboxtr_{P(z)}W_{2}$
is a $V$-module.
By Proposition \ref{pztensor} and
Corollary \ref{pztensor2}, the conclusions of Proposition \ref{pztensor} and
Corollary \ref{pztensor2} are true.

Next we assume in addition that all irreducible $V_{0}$-modules are
${\Bbb R}$-graded and $V$ and $V_{0}$ satisfy Condition 1 and
Condition 2
in Theorem \ref{assoc}. Since any irreducible $V$-module is a $V_{0}$-module
and by assumption any irreducible $V_{0}$-module, thus also any $V_{0}$-module,
is ${\Bbb R}$-graded, we see that any irreducible $V$-module is
${\Bbb R}$-graded. In the proofs of Theorem \ref{assoc} and
Theorem \ref{intwassoc} in \cite{H1} and thus also in
the proof of Theorem \ref{commu} in the preceeding section,
the rationality is used
only in  the proof  of Lemma 14.4 and in Remark 16.1 in \cite{H1}.
So to prove the conclusions of Theorem \ref{assoc},
Theorem \ref{intwassoc} and Theorem \ref{commu} in this case, we need
first to show that the conclusions of Lemma 14.4 and Remark 16.1 in \cite{H1}
 are still true in this case.
It is clear that the conclusions of Lemma 14.4 in \cite{H1} are still
true since ${\cal Y}$ is also an intertwining operators for $V_{0}$ of
type  ${W_{3}}\choose {W_{1}W_{2}}$ if $W_{1}$, $W_{2}$ and $W_{3}$
are viewed as $V_{0}$-modules and thus Lemma 14.4 in \cite{H1} can be
used. Conclusions of Remark 16.1 in \cite{H1} are also true in this case
since
they are obtained by comparing weights and coefficients of series
and thus Remark 16.1
for intertwining operators for $V_{0}$ can be used to obtain
the conclusions in this case.
To show that the conclusions of
Theorem \ref{assoc}, Theorem \ref{intwassoc} and
Theorem \ref{commu} are true in this case,
we still need to show that $V$ satisfies Condition 1 and Condition 2
in  Theorem \ref{assoc}.  Condition 1  is
satisfied by assumption. Since any intertwining operator for $V$ are
also intertwining operators for $V_{0}$ when $V$-modules are viewed as
$V_{0}$-modules, Condition 2 is also satisfied.

Finally we assume in addition that $V_{0}$ satisfies Condition 3 in
Theorem \ref{subvtc}. Since any intertwining operator for $V$ are
also intertwining operators for $V_{0}$ when $V$-modules are viewed as
$V_{0}$-modules, Condition 2 in Theorem \ref{assoc}
and Condition 3 in Theorem \ref{vtc} are
satisfied. Condition 1 in Theorem \ref{assoc} is satisfied by assumption.
We already know that any irreducible $V$-module is ${\Bbb R}$-graded.
As in the proofs of Theorem \ref{assoc}, Theorem \ref{intwassoc} and
Theorem \ref{commu},
in the proof of Theorem \ref{vtc},
the rationality of $V$ is also used only in the proof of
Lemma 14.4 and in Remark 16.1 in \cite{H1}. So the conclusion of
Theorem \ref{vtc} is true in this case.

\renewcommand{\theequation}{\thesection.\arabic{equation}}
\renewcommand{\therema}{\thesection.\arabic{rema}}
\setcounter{equation}{0}
\setcounter{rema}{0}

\section{Proof of Theorem 3.5}

We first prove the case with $m=1$, that is, the case for
minimal Virasoro vertex operator algebras.

Let $W$ be a lower-truncated generalized
$L(c_{p, q}, 0)$-module generated by a homogeneous vector
$w\in W$. Then any element of
$W$ is a linear combination of the elements
\begin{equation}\label{6-1}
L(n_{1})\cdots L(n_{k})w,\;\;\;k\in {\Bbb N}, \;n_{1}, \dots, n_{k}\in
{\Bbb Z}.
\end{equation}
Using the Virasoro commutator relations for the operators
$L(n)$, $n\in {\Bbb Z}$, any element of the form (\ref{6-1}), and thus
any element of $W$,
can be expressed as linear combinations of the elements
\begin{equation}\label{6-2}
L(-m_{1})\cdots L(-m_{k})L(n_{1})\cdots L(n_{l})w,\;\;\;k, l\in {\Bbb N}, \;
m_{1}, \dots, m_{k}, n_{1}, \dots, n_{l}\in {\Bbb Z}_{+}.
\end{equation}
Since $W$ is lower-truncated, we see that for any fixed complex number, there
are only finitely many elements of the form (\ref{6-2}) with weight equal to
this complex number. Thus the homogeneous subspaces of $W$ are all
finite-dimensional, proving that $W$ is a module.
The first conclusion
 is now an
immediate consequence.

We prove the second conclusion  using the differential equations of
Belavin, Polyakov and Zamolodchikov (BPZ equations)
for correlation functions in the
minimal models \cite{BPZ} and the theory of differential equations of
regular singular points. It is clear that we need only to prove the case
that the intertwining operators are among irreducible
$L(c_{p, q}, 0)$-modules.
Given any irreducible
$L(c_{p, q}, 0)$-modules $W_{1}$, $W_{2}$, $W_{3}$, $W_{4}$ and $W_{5}$,
let ${\cal Y}_{1}$, ${\cal Y}_{2}$,
be intertwining operators of type ${W_{4}}\choose {W_{1}W_{5}}$,
${W_{5}}\choose {W_{2}W_{3}}$, respectively and $w_{(1)}$,
$w_{(2)}$, $w_{(3)}$ and $w'_{(4)}$ the lowest weight vectors of
$W_{1}$, $W_{2}$, $W_{3}$ and $W'_{4}$, respectively. We first would like to
show that when $|z_{1}|>|z_{2}|>0$,
\begin{equation}\label{6-3}
\langle w'_{(4)} {\cal Y}_{1}(w_{(1)}, x_{1})
{\cal Y}_{2}(w_{(2)}, x_{2})w_{(3)}\rangle\lbar_{x_{1}^{n}=e^{n\log
z_{1}}, x_{2}^{n}=e^{n \log z_{2}}, n\in \Bbb{C}}
\end{equation}
 satisfies a BPZ equation. {}From the representation theory of the Virasoro
algebra, we know that there exists
\begin{equation}\label{6-3.5}
P=\sum_{i=1}^{k}a_{i}L(-m^{(i)}_{1})\cdots L(-m^{(i)}_{l_{i}})\in
U({\frak L}_{-})
\end{equation}
where $l_{i}, m^{(i)}_{j}\in {\Bbb Z}_{+}$, $1\le j\le l_{i}$,
$1\le i\le k$ such that
\begin{equation}\label{6-3.6}
\sum_{j=1}^{l_{1}}m^{(1)}_{j}=\cdots =\sum_{j=1}^{l_{k}}m^{(k)}_{j}>0.
\end{equation}
and
$$
Pw_{(3)}=0.
$$
Thus we have
\begin{equation}\label{6-4}
\langle w'_{(4)} {\cal Y}_{1}(w_{(1)}, x_{1})
{\cal Y}_{2}(w_{(2)}, x_{2})Pw_{(3)}\rangle\lbar_{x_{1}^{n}=e^{n\log
z_{1}}, x_{2}^{n}=e^{n \log z_{2}}, n\in \Bbb{C}}=0.
\end{equation}

Since $w_{(1)}$, $w_{(2)}$ and $w'_{(4)}$ are lowest weight vectors,
\begin{eqnarray}\label{6-5}
\lefteqn{\langle w'_{(4)} {\cal Y}_{1}(w_{(1)}, x_{1})
{\cal Y}_{2}(w_{(2)}, x_{2})L(-m)w\rangle\lbar_{x_{1}^{n}=e^{n\log
z_{1}}, x_{2}^{n}=e^{n \log z_{2}}, n\in \Bbb{C}}=}\nno\\
&&=-\langle w'_{(4)} (x_{1}^{-m+1}\frac{\p}{\p x_{1}}
+(\wt w_{(1)})(-m+1)x_{1}^{-m})\cdot\nno\\
&&\quad\quad\cdot {\cal Y}_{1}(w_{(1)}, x_{1})
{\cal Y}_{2}(w_{(2)}, x_{2})w\rangle\lbar_{x_{1}^{n}=e^{n\log
z_{1}}, x_{2}^{n}=e^{n \log z_{2}}, n\in \Bbb{C}}\nno\\
&&\quad -\langle w'_{(4)}
{\cal Y}_{1}(w_{(1)}, x_{1})
 (x_{2}^{-m+1}\frac{\p}{\p x_{2}}+(\wt w_{(2)})(-m+1)x_{2}^{-m})\cdot\nno\\
&&\quad\quad\cdot {\cal Y}_{2}(w_{(2)}, x_{2})w\rangle\lbar_{x_{1}^{n}=e^{n\log
z_{1}}, x_{2}^{n}=e^{n \log z_{2}}, n\in \Bbb{C}}\nno\\
&&=-(z_{1}^{-m+1}\frac{\p}{\p z_{1}}
+(\wt w_{(1)})(-m+1)z_{1}^{-m}\nno\\
&&\quad\quad
+(z_{2}^{-m+1}\frac{\p}{\p z_{2}}+(\wt w_{(2)})(-m+1)z_{2}^{-m}))\cdot\nno\\
&&\quad\quad\quad\cdot \langle w'_{(4)}
{\cal Y}_{1}(w_{(1)}, x_{1})
{\cal Y}_{2}(w_{(2)}, x_{2})w\rangle\lbar_{x_{1}^{n}=e^{n\log
z_{1}}, x_{2}^{n}=e^{n \log z_{2}}, n\in \Bbb{C}}
\end{eqnarray}
for any $w\in W_{(3)}$ and $m \in {\Bbb Z}_{+}$.
{}From (\ref{6-3.5}), (\ref{6-4}) and (\ref{6-5}), we see that (\ref{6-3})
satisfies a differential equation.
Similarly there exists
\begin{equation}\label{6-6}
Q=\sum_{i=1}^{r}b_{i}L(-n^{(i)}_{1})\cdots L(-n^{(i)}_{s_{i}})\in
U({\frak L}_{-})
\end{equation}
where $s_{i}, n^{(i)}_{j}\in {\Bbb Z}_{+}$, $1\le j\le s_{i}$,
$1\le i\le r$ such that
\begin{equation}\label{6-7}
\sum_{j=1}^{s_{1}}n^{(1)}_{j}=\cdots =\sum_{j=1}^{s_{r}}n^{(r)}_{j}>0.
\end{equation}
and
$$
Qw_{(2)}=0.
$$
Thus we have
\begin{equation}\label{6-9}
\langle w'_{(4)} {\cal Y}_{1}(w_{(1)}, x_{1})
{\cal Y}_{2}(Qw_{(2)}, x_{2})w_{(3)}\rangle\lbar_{x_{1}^{n}=e^{n\log
z_{1}}, x_{2}^{n}=e^{n \log z_{2}}, n\in \Bbb{C}}=0.
\end{equation}

Recalling that we have a linear isomorphism
$$\Omega_{-1}: {\cal V}_{W_{2}W_{3}}^{W_{5}}\to {\cal V}_{W_{3}W_{2}}^{W_{5}}$$
and its inverse
$$\Omega_{0}: {\cal V}_{W_{3}W_{2}}^{W_{5}}\to {\cal V}_{W_{2}W_{3}}^{W_{5}}$$
defined by
$$\Omega_{0}({\cal Y}(w_{(3)}, x)w_{(2)}=e^{xL(-1)}{\cal Y}(w_{(2)},
e^{\pi i}x)w_{(3)}$$
(see \cite{FHL} and \cite{HL3}), we obtain
\begin{eqnarray}\label{6-10}
\lefteqn{\langle w'_{(4)}, {\cal Y}_{1}(w_{(1)}, x_{1})
{\cal Y}_{2}(Qw_{(2)}, x_{2})w_{(3)}\rangle\lbar_{x_{1}^{n}=e^{n\log
z_{1}}, x_{2}^{n}=e^{n \log z_{2}}, n\in \Bbb{C}}=}\nno\\
&&=\langle w'_{(4)}, {\cal Y}_{1}(w_{(1)}, x_{1})
\Omega_{0}(\Omega_{-1}({\cal Y}_{2}))(Qw_{(2)}, x_{2})w_{(3)}
\rangle\lbar_{x_{1}^{n}=e^{n\log
z_{1}}, x_{2}^{n}=e^{n \log z_{2}}, n\in \Bbb{C}}\nno\\
&&=\langle w'_{(4)}, {\cal Y}_{1}(w_{(1)}, x_{1})e^{x_{2}L(-1)}\cdot\nno\\
&&\quad\quad\cdot \Omega_{-1}({\cal Y}_{2})(w_{(3)}, e^{\pi i}x_{2})Qw_{(2)}
\rangle\lbar_{x_{1}^{n}=e^{n\log
z_{1}}, x_{2}^{n}=e^{n \log z_{2}}, n\in \Bbb{C}}\nno\\
&&=\langle e^{x_{2}L(1)}w'_{(4)}, {\cal Y}_{1}(w_{(1)}, x_{1}-x_{2})\cdot\nno\\
&&~~~~~~
\cdot\Omega_{-1}~({\cal Y}_{2})(w_{(3)},
e^{\pi i}x_{2})Qw_{(2)}
\rangle\lbar_{x_{1}^{n}=e^{n\log
z_{1}}, x_{2}^{n}=e^{n \log z_{2}}, n\in \Bbb{C}}\nno\\
&&=\langle w'_{(4)}, {\cal Y}_{1}(w_{(1)}, x_{1}-x_{2})\nno\\
&&\quad\quad\cdot \Omega_{-1}({\cal Y}_{2})(w_{(3)}, e^{\pi i}x_{2})
Qw_{(2)}
\rangle\lbar_{x_{1}^{n}=e^{n\log
z_{1}}, x_{2}^{n}=e^{n \log z_{2}}, n\in \Bbb{C}}.
\end{eqnarray}

Using (\ref{6-6}), (\ref{6-9}) and (\ref{6-10}) and the same method
as above, we obtain another differential
equation satisfied by (\ref{6-3}).
{}From (\ref{6-5}) and (\ref{6-10}),
we see that  $z_{1}=z_{2}$ are not singularities of the
differential equation obtained {}from
(\ref{6-4}) but are singularities of the differential equation obtained {}from
(\ref{6-9}). So the two differential equations are independent.
We obtain a system of two independent differential equations for
functions with two variables. Formally, (\ref{6-3}) satisfies this
system of equations. {}From
(\ref{6-3.6}), (\ref{6-5}), (\ref{6-7}) and (\ref{6-10}), we see that
this system has only regular singularities $z_{1}, z_{2}=0, \infty$
and $z_{1}=z_{2}$. Since ${\cal Y}_{1}$ and ${\cal Y}_{2}$ are intertwining
operators among irreducible modules, there exist rational numbers
$h_{1}$ and $h_{2}$ such that
$${\cal Y}_{1}(w_{(1)}, x_{1})\in x^{h_{1}}\mbox{\rm Hom}(W_{5},
W_{4})[[x_{1}, x_{1}^{-1}]]$$ and
$${\cal Y}_{2}(w_{(2)}, x_{2})\in x^{h_{2}}\mbox{\rm Hom}(W_{3},
W_{5})[[x_{2}, x_{2}^{-1}]].$$ Thus by the definition of intertwining
operator, there are two rational numbers $t_{1}$, $t_{2}$ and a
power series $g(z)$ such that
(\ref{6-3}) is equal to
\begin{equation}\label{6-11}
z_{1}^{t_{1}}z_{2}^{t_{2}}g(z_{2}/z_{1}).
\end{equation}
The system of differential equations for (\ref{6-3}) give a
differential equation
for $g(z)$. Since the system of
differential equations for (\ref{6-3}) only have
singularity when $z_{1}, z_{2}=0, \infty$ and $z_{1}=z_{2}$ and
$g(z)$ is a power series, this equation
for $g(z)$ has coefficients analytic in $|z|<1$.
The coefficients of $g(z)$ are can be calculated using ${\cal Y}_{1}$,
${\cal Y}_{2}$, $w_{(1)}$, $w_{(2)}$, $w_{(3)}$ and $w'_{(4)}$ and
thus the initial conditions at $z=0$ for the initial-value problem for
the differential equation for $g(z)$
are known. By the existence
and uniqueness theorem of differential equations with analytic coefficients,
$g(z)$ must be convergent when $|z|<1$. Thus (\ref{6-11}) and (\ref{6-3}) are
convergent when $|z_{2}|>|z_{1}|>0$.

Now we can use the values of
(\ref{6-3}) and its derivatives
at $(z_{1}, z_{2})$ satisfying $|z_{2}|>|z_{1}|>|z_{1}-z_{2}|>0$
as the initial condition to solve the system of equations
of regular singular points for (\ref{6-3}) in the region
$|z_{2}|>|z_{1}-z_{2}|>0$.
By the theory of equations of regular
singular points (see, for example, Appendix B of \cite{K}),
there exist
$j\in \Bbb{N}$, $r_{i}, s_{i}\in \Bbb{C}$, $i=1, \dots, j$,
$M\in {\Bbb N}$ and analytic
functions $f_{i, m_{1}, m_{2}}(z)$ on $|z|<1$, $i=1, \dots, j$,
$m_{1}, m_{2}\in {\Bbb N}$, $m_{1}+m_{2}\le M$, such that
the solution in
this region
is of the form
\begin{eqnarray}\label{6-12}
\lefteqn{\sum_{i=1}^{j}\sum_{m_{1}, m_{2}\in {\Bbb N}, m_{1}+m_{2}\le M}
z_{2}^{r_{i}}(z_{1}-z_{2})^{s_{i}}
(\log z_{2})^{m_{1}}}\cdot\nno\\
&&\quad\quad\quad\quad\quad\quad\cdot
\left(\log \left(\frac{z_{1}-z_{2}}{z_{2}}\right)\right)^{m_{2}}
f_{i, m_{1}, m_{2}}\left(\frac{z_{1}-z_{2}}{z_{2}}\right).
\end{eqnarray}
Since (\ref{6-12}) is the analytic extension of (\ref{6-3}) and thus also
of (\ref{6-11}) and  since
(\ref{6-11}) at the  singularity $z_{2}=0$ does not contain
any term involving $\log z_{2}$, we have
$$f_{i, m_{1}, m_{2}}(z)=0$$
for any $i=1, \dots, j$ and any $z$ satisfying
$|z|<1$ if $m_{1}$ or $m_{2}$ is not $0$.
Comparing (\ref{6-12}) with (\ref{6-11}), we also see that we can take
$r_{i}$ and $s_{i}$, $i=1, \dots, j$, to be rational numbers.
Thus (\ref{6-12}) or the
analytic extension of (\ref{6-3}) in the region $|z_{2}|>|z_{1}-z_{2}|>0$
must be of the form (\ref{phyper}).

Let $N$ be an integer such that
for this extension and $w_{(1)}$ and $w_{(2)}$ above,
(\ref{si}) holds.
We now use induction on the weights of $w_{(1)}$, $w_{(2)}$, $w_{(3)}$ and
$w'_{(4)}$ to show that for any $w_{(1)}$, $w_{(2)}$, $w_{(3)}$ and
$w'_{(4)}$, not necessarily lowest weight vectors, (\ref{6-3}) is
convergent when $|z_{1}|>|z_{2}|>0$ and can be analytically extended to
a function of the form (\ref{phyper}) in the region $|z_{2}|>|z_{1}-z_{2}|>0$
satisfying (\ref{si}) for the $N$ above.
First we assume that $w_{(1)}$ and $w_{(2)}$ are still lowest weight vectors.
In this case, instead of (\ref{6-5}), we have
\begin{eqnarray}\label{6-14}
\lefteqn{\langle w'_{(4)}, {\cal Y}_{1}(w_{(1)}, x_{1})
{\cal Y}_{2}(w_{(2)}, x_{2})L(-m)w\rangle\lbar_{x_{1}^{n}=e^{n\log
z_{1}}, x_{2}^{n}=e^{n \log z_{2}}, n\in \Bbb{C}}}\nno\\
&&=\langle L(m)w'_{(4)},
{\cal Y}_{1}(w_{(1)}, x_{1})
{\cal Y}_{2}(w_{(2)}, x_{2})w\rangle\lbar_{x_{1}^{n}=e^{n\log
z_{1}}, x_{2}^{n}=e^{n \log z_{2}}, n\in \Bbb{C}}\nno\\
&& \quad -(z_{1}^{-m+1}\frac{\p}{\p z_{1}}
+(\wt w_{(1)})(-m+1)z_{1}^{-m}\nno\\
&&\quad\quad  +z_{2}^{-m+1}\frac{\p}{\p z_{2}}+(\wt w_{(2)})(-m+1)z_{2}^{-m})
\cdot\nno\\
&&\quad\quad \quad\cdot \langle w'_{(4)},
{\cal Y}_{1}(w_{(1)}, x_{1})
{\cal Y}_{2}(w_{(2)}, x_{2})w\rangle\lbar_{x_{1}^{n}=e^{n\log
z_{1}}, x_{2}^{n}=e^{n \log z_{2}}, n\in \Bbb{C}}.\hspace{2em}
\end{eqnarray}
Since the weight of $L(m)w'_{(4)}$ is less than the weight of $w'_{(4)}$ and
the weight of $w$ is less than the weight of $L(-m)w$, using induction
we see that the left-hand side of (\ref{6-14}) is
convergent when $|z_{1}|>|z_{2}|>0$ and can be analytically extended to
a function of the form (\ref{phyper}) in the region $|z_{2}|>|z_{1}-z_{2}|>0$.
Also {}from (\ref{6-14}), we see that if
 the analytic extensions of
$$\langle L(m)w'_{(4)},
{\cal Y}_{1}(w_{(1)}, x_{1})
{\cal Y}_{2}(w_{(2)}, x_{2})w\rangle\lbar_{x_{1}^{n}=e^{n\log
z_{1}}, x_{2}^{n}=e^{n \log z_{2}}, n\in \Bbb{C}}$$
and
$$\langle w'_{(4)},
{\cal Y}_{1}(w_{(1)}, x_{1})
{\cal Y}_{2}(w_{(2)}, x_{2})w\rangle\lbar_{x_{1}^{n}=e^{n\log
z_{1}}, x_{2}^{n}=e^{n \log z_{2}}, n\in \Bbb{C}}$$
 in the region $|z_{2}|>|z_{1}-z_{2}|>0$
are of the form (\ref{phyper}) satisfying (\ref{si}),
the extension of the left-hand side of (\ref{6-14}) can also
be written in the form (\ref{phyper}) satisfying (\ref{si}). Similarly we
can prove the case for $w'_{(4)}=L'(-m)w'$ using the case for $w'$.

Next we consider the case that only $w_{(1)}$ is still
a lowest weight vector. By the Jacobi identity and the definition
of contragredient vertex operator, we have
\begin{eqnarray}\label{6-15}
\lefteqn{\langle w'_{(4)},
{\cal Y}_{1}(w_{(1)}, x_{1})
{\cal Y}_{2}(L(-m)w, x_{2})w_{(3)}\rangle\lbar_{x_{1}^{n}=e^{n\log
z_{1}}, x_{2}^{n}=e^{n \log z_{2}}, n\in \Bbb{C}}}\nno\\
&&=-\res_{x}(-x+x_{2})^{-m+1}\langle w'_{(4)},
{\cal Y}_{1}(w_{(1)}, x_{1})\cdot\nno\\
&&\quad\quad\cdot
{\cal Y}_{2}(w, x)Y_{3}(\omega, x_{2})w_{(3)}\rangle\lbar_{x_{1}^{n}=e^{n\log
z_{1}}, x_{2}^{n}=e^{n \log z_{2}}, n\in \Bbb{C}}\nno\\
&&\quad +\res_{x}(x_{2}-x)^{-m+1}\langle w'_{(4)},
Y_{4}(\omega, x)\cdot\nno\\
&&\quad\quad\cdot{\cal Y}_{1}(w_{(1)}, x_{1})
{\cal Y}_{2}(w, x)w_{(3)}\rangle\lbar_{x_{1}^{n}=e^{n\log
z_{1}}, x_{2}^{n}=e^{n \log z_{2}}, n\in \Bbb{C}}\nno\\
&&\quad -\res_{x}(x-x_{2})^{-m+1}
\res_{x_{0}}x^{-1}_{0}\delta\left(\frac{x-x_{0}}{x_{1}}\right)
\langle w'_{(4)},
{\cal Y}_{1}(Y_{1}(\omega, x_{0})w_{(1)}, x_{1})\cdot\nno\\
&&\quad\quad\cdot
{\cal Y}_{2}(w, x_{2})w_{(3)}\rangle\lbar_{x_{1}^{n}=e^{n\log
z_{1}}, x_{2}^{n}=e^{n \log z_{2}}, n\in \Bbb{C}}\nno\\
&&=-\res_{x}(-x_{2}+x)^{-m+1}\langle w'_{(4)},
{\cal Y}_{1}(w_{(1)}, x_{1})\cdot\nno\\
&&\quad\quad\cdot
{\cal Y}_{2}(w, x_{2})Y_{3}(\omega, x)w_{(3)}\rangle\lbar_{x_{1}^{n}=e^{n\log
z_{1}}, x_{2}^{n}=e^{n \log z_{2}}, n\in \Bbb{C}}\nno\\
&&\quad +\res_{x}(x-x_{2})^{-m+1}x^{-2}\langle Y'_{4}(\omega, x^{-1})w'_{(4)},
{\cal Y}_{1}(w_{(1)}, x_{1})\cdot\nno\\
&&\quad\quad\cdot
{\cal Y}_{2}(w, x)w_{(3)}\rangle\lbar_{x_{1}^{n}=e^{n\log
z_{1}}, x_{2}^{n}=e^{n \log z_{2}}, n\in \Bbb{C}}\nno\\
&&\quad -\res_{x}(x-x_{2})^{-m+1}
\res_{x_{0}}x^{-1}_{0}\delta\left(\frac{x-x_{0}}{x_{1}}\right)
\langle w'_{(4)},
{\cal Y}_{1}((x^{-2}_{0}L(0)\nno\\
&&\quad\quad +x_{0}^{-1}L(-1))w_{(1)}, x_{1})
{\cal Y}_{2}(w, x_{2})w_{(3)}\rangle\lbar_{x_{1}^{n}=e^{n\log
z_{1}}, x_{2}^{n}=e^{n \log z_{2}}, n\in \Bbb{C}}\nno\\
&&=-\res_{x}(-x_{2}+x)^{-m+1}\langle w'_{(4)},
{\cal Y}_{1}(w_{(1)}, x_{1})\cdot\nno\\
&&\quad\quad\cdot
{\cal Y}_{2}(w, x_{2})Y_{3}(\omega, x)w_{(3)}\rangle\lbar_{x_{1}^{n}=e^{n\log
z_{1}}, x_{2}^{n}=e^{n \log z_{2}}, n\in \Bbb{C}}\nno\\
&&\quad +\res_{x}(x-x_{2})^{-m+1}x^{-2}\langle Y'_{4}(\omega, x^{-1})w'_{(4)},
{\cal Y}_{1}(w_{(1)}, x_{1})\cdot\nno\\
&&\quad\quad\cdot
{\cal Y}_{2}(w, x)w_{(3)}\rangle\lbar_{x_{1}^{n}=e^{n\log
z_{1}}, x_{2}^{n}=e^{n \log z_{2}}, n\in \Bbb{C}}\nno\\
&&\quad -(\wt w_{(1)})\res_{x}(x-x_{2})^{-m+1}
\res_{x_{0}}x^{-3}_{0}\delta\left(\frac{x-x_{0}}{x_{1}}\right)
\langle w'_{(4)},
{\cal Y}_{1}(w_{(1)}, x_{1})\cdot\nno\\
&&\quad\quad\cdot
{\cal Y}_{2}(w, x_{2})w_{(3)}\rangle\lbar_{x_{1}^{n}=e^{n\log
z_{1}}, x_{2}^{n}=e^{n \log z_{2}}, n\in \Bbb{C}}\nno\\
&&\quad -\res_{x}(x-x_{2})^{-m+1}
\res_{x_{0}}x^{-2}_{0}\delta\left(\frac{x-x_{0}}{x_{1}}\right)
\frac{\p}{\p x_{1}}\langle w'_{(4)},
{\cal Y}_{1}(w_{(1)}, x_{1})\cdot\nno\\
&&\quad\quad\cdot
{\cal Y}_{2}(w, x_{2})w_{(3)}\rangle\lbar_{x_{1}^{n}=e^{n\log
z_{1}}, x_{2}^{n}=e^{n \log z_{2}}, n\in \Bbb{C}}
\end{eqnarray}
for any $w\in W_{2}$ and $m\in {\Bbb Z}_{+}$. Using induction on the
weight of $w_{(2)}$, we see {}from
(\ref{6-15}) that the left-hand side of of (\ref{6-15}) is
convergent when $|z_{1}|>|z_{2}|>0$ and can be analytically extended to
a function of the form (\ref{phyper}) in the region $|z_{2}|>|z_{1}-z_{2}|>0$
satisfying (\ref{si}). The most general case that $w_{(1)}$, $w_{(2)}$,
$w_{(3)}$, $w'_{(4)}$ are arbitrary can be proved similarly using
a formula a little more complicated than (\ref{6-1}) and induction
on the weight of $w_{(1)}$.

The proof of the third conclusion in this special case
 is similar to the proof of the
convergence in the proof of the second conclusion above  by deriving certain
BPZ equations.

To prove the general case,  we need:

\begin{lemma}\label{intwop}
Let $m$ be a positive integer, $(p_{i}, q_{i})$, $i=1, \dots, m$, $m$
pairs of relatively prime positive integers larger than $1$,
$V=L(c_{p_{1}, q_{1}}, 0)\otimes \cdots \otimes L(c_{p_{m}, q_{m}},
0)$, $W_{i}=L(c_{p_{1}, q_{1}}, h_{1}^{(t)})\otimes \cdots
\otimes L(c_{p_{m}, q_{m}}, h_{m}^{(t)})$,
$t=1, 2, 3$, irreducible $V$-modules and ${\cal Y}$ an
intertwining operator of type ${W_{3}}\choose {W_{1}W_{2}}$. Then
there exist intertwining operators ${\cal Y}_{i}$
of type ${L(c_{p_{i}, q_{i}}, h_{i}^{(3)})} \choose
{L(c_{p_{i}, q_{i}}, h_{i}^{(1)})L(c_{p_{i}, q_{i}}, h_{i}^{(2)})}$,
$i=1, \dots, m$, such that
\begin{eqnarray}\label{6-16}
{\cal Y}={\cal Y}_{1}\otimes \cdots \otimes {\cal Y}_{m}.
\end{eqnarray}
\end{lemma}
\pf
Since $L(c_{p_{i}, q_{i}}, 0)$, $i=1, \dots, m$, are rational, by
Proposition 2.10 in \cite{DMZ}, we know that the vector space
${\cal V}_{W_{1}W_{2}}^{W_{3}}$ is isomorphic the vector space
\begin{equation}\label{6-17}
\otimes_{i=1}^{m}
{\cal V}_{L(c_{p_{i}, q_{i}}, h_{i}^{(1)})L(c_{p_{i}, q_{i}}, h_{i}^{(2)})}
^{L(c_{p_{i}, q_{i}}, h_{i}^{(3)})}.
\end{equation}
On the other hand, intertwining operators of the form (\ref{6-16})
form a subspace of ${\cal V}_{W_{1}W_{2}}^{W_{3}}$ with dimension
equal to the dimension of (\ref{6-17}). So this subspace must equal to
${\cal V}_{W_{1}W_{2}}^{W_{3}}$ itself, proving the lemma.
This lemma can also be proved directly
using the special properties of the Virasoro vertex operator algebras.
\epfv

Now we prove the theorem in the general case.

Since $L(c_{p_{i}, q_{i}}, 0)$, $i=1, \dots, m$, are rational,
$V=L(c_{p_{1}, q_{1}}, 0)
\otimes \cdots \otimes
L(c_{p_{n}, q_{m}}, 0)$
is also rational \cite{DMZ}. Since for any $i\in {\Bbb Z}$, $1\le i\le m$,
every finitely-generated
lower-truncated generalized
$L(c_{p_{i}, q_{i}}, 0)$-module is a module,
the same method proving the rationality
of $V=L(c_{p_{1}, q_{1}}, 0)
\otimes \cdots \otimes
L(c_{p_{m}, q_{m}}, 0)$ shows that any finitely-generated
lower-truncated generalized $V$-module $W$ is a
sum of $V$-modules. Since there are only finitely
many irreducible $V$-modules and $W$
is finitely-generated, $W$ must be an $V$-module.

For simplicity we only prove the second conclusion for $m=2$. Consider
intertwining operators ${\cal Y}_{1}$ and ${\cal Y}_{2}$ of type
${W_{4}}\choose {W_{1}W_{5}}$ and ${W_{5}}\choose {W_{2}W_{3}}$, respectively.
For simplicity, we assume that $W_{1}, \dots, W_{5}$ are irreducible
modules. Thus we have $W_{i}=L(c_{p_{1}, q_{1}}, h_{1}^{(t)})\otimes
L(c_{p_{2}, q_{2}}, h_{2}^{(t)})$, $t=1, \dots, 5$. By Lemma \ref{intwop},
there are intertwining operators
${\cal Y}_{1}^{(i)}$ and ${\cal Y}_{1}^{(i)}$, $i=1, 2$,
of type ${L(c_{p_{t}, q_{i}}, h_{i}^{(5)})}
\choose {L(c_{p_{i}, q_{i}}, h_{i}^{(1)})L(c_{p_{i}, q_{i}}, h_{i}^{(4)})}$ and
${L(c_{p_{i}, q_{i}}, h_{i}^{(4)})}
\choose {L(c_{p_{i}, q_{i}}, h_{i}^{(2)})L(c_{p_{i}, q_{i}}, h_{i}^{(3)})}$,
respectively,
such that ${\cal Y}_{j}={\cal Y}_{j}^{(1)}\otimes {\cal Y}_{j}^{(2)}$,
$j=1, 2$. By the case with $m=1$, the products of the intertwining
operators for
$L(c_{p_{i}, q_{i}}, 0)$, $i=1, 2$,
 have the convergence and extension property.
Let $N_{1}$ and $N_{2}$ be the integers
for ${\cal Y}_{1}^{(1)}, {\cal Y}_{2}^{(1)}$ and
${\cal Y}_{1}^{(2)}, {\cal Y}_{2}^{(2)}$, respectively,
in the convergence and extension
property for products. Then for any
$w^{(i)}_{(l)}\in L(c_{p_{i}, q_{i}}, h_{i}^{(t)})$, $i=1, 2$, $t=1, 2, 3$ and
$(w_{(4)}^{(i)})'\in L(c_{p_{i}, q_{i}}, h_{i}^{(4)})'$, $i=1, 2$, there exist
rational numbers $r^{(i)}_{j}, s^{(i)}_{j}$ and analytic
functions $f^{(i)}_{j}(z)$ on $|z|<1$, $j=1, \dots, k_{i}$, $i=1, 2$,
satisfying
$\wt w^{(i)}_{(1)}+\wt w_{(2)}^{(i)}+s^{(i)}_{j}>N_{i}$ such that
$$
\langle (w^{(i)}_{(4)})',
{\cal Y}_{1}^{(i)}(w^{(i)}_{(1)}, x_{1})
({\cal Y}^{(i)}_{2}(w^{(i)}_{(2)}, x_{2}) w^{(i)}_{(3)}
\rangle_{L(c_{p_{(i)}, q_{(i)}}, h_{(i)}^{(4)})}
\lbar_{x_{1}^{n}=e^{n\log z_{1}},\; x_{2}^{n}=e^{n\log z_{2}},
\;n\in {\Bbb Q}}
$$
is absolutely convergent when $|z_{1}|>|z_{2}|>0$
 and can be analytically extended to
$$
\sum_{j=1}^{k}z_{2}^{r^{(i)}_{j}}(z_{1}-z_{2})^{s^{(i)}_{j}}
f^{(i)}_{j}\left(\frac{z_{1}-z_{2}}{z_{2}}\right)
$$
when $|z_{2}|>|z_{1}-z_{2}|>0$.
Thus
\begin{eqnarray*}
\lefteqn{\langle (w^{(1)}_{(4)})'\otimes (w^{(2)}_{(4)})',
{\cal Y}_{(1)}(w^{(1)}_{(1)}\otimes w^{(2)}_{(1)}, x_{1})\cdot}\nno\\
&&\quad\quad\cdot{\cal Y}_{(2)}(w^{(1)}_{(2)}\otimes w^{(2)}_{(2)}, x_{2})
(w^{(1)}_{(3)}w^{(2)}_{(3)})
\rangle_{W_{4}}\lbar_{x_{1}^{n}=e^{n\log z_{1}},\; x_{2}^{n}=e^{n\log z_{2}},
\;n\in {\Bbb Q}}\\
&&=\langle (w^{(1)}_{(4)})',
{\cal Y}_{1}^{(1)}(w^{(1)}_{(1)}, x_{1})\cdot\nno\\
&&\quad\quad\cdot
{\cal Y}^{(1)}_{2}(w^{(1)}_{(2)}, x_{2}) w^{(1)}_{(3)}
\rangle_{L(c_{p_{1}, q_{1}}, h_{1}^{(4)})}
\lbar_{x_{1}^{n}=e^{n\log z_{1}},\; x_{2}^{n}=e^{n\log z_{2}},
\;n\in {\Bbb Q}}\\
&&\quad\quad\cdot \langle (w^{(2)}_{(4)})',
{\cal Y}_{(1)}^{2}(w^{(2)}_{(1)}, x_{1})\cdot\nno\\
&&\quad\quad\cdot
{\cal Y}^{(2)}_{2}(w^{(2)}_{(2)}, x_{2}) w^{(2)}_{(3)}
\rangle_{L(c_{p_{2}, q_{2}}, h_{2}^{(4)})}
\lbar_{x_{1}^{n}=e^{n\log z_{1}},\; x_{2}^{n}=e^{n\log z_{2}},
\;n\in {\Bbb Q}}
\end{eqnarray*}
is absolutely convergent  and can be analytic extended to
\begin{eqnarray*}
\lefteqn{\sum_{j_{1}=1}^{k_{1}}z_{2}^{r^{(1)}_{j_{1}}}
(z_{1}-z_{2})^{s^{(1)}_{j_{1}}}
f^{(1)}_{j_{1}}\left(\frac{z_{1}-z_{2}}{z_{2}}\right)\cdot }\nno\\
&&\quad\cdot\sum_{j_{2}=1}^{k_{2}}z_{2}^{r^{(2)}_{j_{2}}}
(z_{1}-z_{2})^{s^{(2)}_{j_{2}}}
f^{(2)}_{j_{2}}\left(\frac{z_{1}-z_{2}}{z_{2}}\right)\\
&&=\sum_{j_{1}=1}^{k_{1}}\sum_{j_{2}=1}^{k_{2}}
z_{2}^{(r^{(1)}_{j_{1}}+r^{(2)}_{j_{2}})}
(z_{1}-z_{2})^{(s^{(1)}_{j_{1}}+s^{(2)}_{j_{2}})}
f^{(1)}_{j_{1}}\left(\frac{z_{1}-z_{2}}{z_{2}}\right)
f^{(2)}_{j_{2}}\left(\frac{z_{1}-z_{2}}{z_{2}}\right).
\end{eqnarray*}
And we have
\begin{eqnarray*}
\lefteqn{\wt (w^{(1)}_{(1)}\otimes w^{(2)}_{(1)})
+\wt (w^{(1)}_{(2)}\otimes w^{(2)}_{(2)})+(s^{(1)}_{j_{1}}+s^{(2)}_{j_{2}})}\\
&&=(\wt w^{(1)}_{(1)}+\wt w_{(2)}^{(1)}+s^{(1)}_{j_{1}})
+(\wt w^{(2)}_{(1)}+\wt w_{(2)}^{(2)}+s^{(2)}_{j_{2}})\\
&&>N_{1}+N_{2},
\end{eqnarray*}
$j_{1}=1, \dots, k_{1}$, $j_{2}=1, \dots, k_{2}$.

 The proof of the third
 conclusion is similar to the proof of the
convergence in the proof of the
second conclusion above.

\renewcommand{\theequation}{\thesection.\arabic{equation}}
\renewcommand{\therema}{\thesection.\arabic{rema}}
\setcounter{equation}{0}
\setcounter{rema}{0}

\section{Proof of Proposition 3.7}

We still only prove the result for $m=2$.
Since $V$ has a subalgebra isomorphic to $L(c_{p_{1}, q_{1}}, 0)\otimes
L(c_{p_{2}, q_{2}}, 0)$,
any  $V$-module is a
$L(c_{p_{1}, q_{1}}, 0)\otimes
L(c_{p_{2}, q_{2}}, 0)$-module. In particular, $V$ is an
$L(c_{p_{1}, q_{1}}, 0)\otimes
L(c_{p_{2}, q_{2}}, 0)$-module. Since $L(c_{p_{1}, q_{1}}, 0)\otimes
L(c_{p_{2}, q_{2}}, 0)$ is
rational, $V$ can be decomposed as a direct sum of irreducible
$L(c_{p_{1}, q_{1}}, 0)\otimes
L(c_{p_{2}, q_{2}}, 0)$-modules. Thus by Theorem 4.7.4 of \cite{FHL}, as a
$L(c_{p_{1}, q_{1}}, 0)\otimes
L(c_{p_{2}, q_{2}}, 0)$-module
$$V=\sum_{j=1}^{k}L(c_{p_{1}, q_{1}}, h_{1}^{(j)})
\otimes L(c_{p_{2}, q_{2}}, h_{2}^{(j)}).$$
Let $W$ be a
lower-truncated generalized $V$-module generated by one homogeneous
element $w\in W$.
Then it is a generalized
$L(c_{p_{1}, q_{1}}, 0)\otimes
L(c_{p_{2}, q_{2}}, 0)$-module.
By a lemma of Dong-Mason \cite{DM} and Li
\cite{L}, $W$ is spanned by elements of the form
$v_{n}w$ where $v\in V$, $n\in {\Bbb C}$.  Let
$u_{(j)}^{(i)}$, $i=1, 2$, $j=1, \dots, k$, be the lowest weight vectors
of $L(c_{p, q}, h_{i}^{(j)})$, respectively. Using the Jacobi identity
 for
generalized modules, we see that elements of the form $v_{n}w$ are
spanned by elements of the form
\begin{eqnarray}
\lefteqn{(L(-m^{(1)}_{1})\cdots L(-m^{(1)}_{p_{1}})\otimes L(-m^{(2)}_{1})
\cdots L(-m^{(2)}_{p_{2}}))\cdot}\nno\\
&&\quad \cdot(L(-1)^{l_{1}}u_{(j)}^{(1)})_{j_{1}}\otimes
(L(-1)^{l_{2}}u_{(j)}^{(2)})_{j_{2}}\cdot \nno\\
&&\quad\quad  \cdot(L(n^{(1)}_{1})\cdots L(n^{(1)}_{q_{1}})\otimes
L(n^{(2)}_{1})
\cdots L(n^{(2)}_{q_{2}}))w,\label{7-1}
\end{eqnarray}
$m^{(1)}_{1}, \dots, m^{(1)}_{p_{1}}$, $m^{(2)}_{1}, \dots, m^{(2)}_{p_{2}}$,
$n^{(1)}_{1}, \dots, n^{(1)}_{p_{1}}$, $n^{(2)}_{1}, \dots, n^{(2)}_{p_{2}}
\in {\Bbb Z}_{+}$,
$l_{1}, l_{2}\in {\Bbb N}$, $j_{1}, j_{2}\in {\Bbb Q}$,
$j=1, \dots, k$. By Theorem
\ref{cpqm},
the generalized $L(c_{p_{1}, q_{1}}, 0)\otimes
L(c_{p_{2}, q_{2}}, 0)$-module generated by $w$ is a
module and thus a finite direct sum of irreducible
$L(c_{p_{1}, q_{1}}, 0)\otimes
L(c_{p_{2}, q_{2}}, 0)$-modules. Thus by Theorem 4.7.4 of \cite{FHL},
$w=\sum_{t=1}^{r}w_{(t)}^{(1)}\otimes
w_{(t)}^{(2)}$ where $w_{(t)}^{(i)}$, $t=1, \dots, r$, $i=1, 2$, are
homogeneous
elements of irreducible $L(c_{p_{i}, q_{i}}, 0)$-modules. Thus elements
of the form (\ref{7-1}) are spanned by
elements of the form
\begin{eqnarray*}
\lefteqn{L(-m^{(1)}_{1})\cdots L(-m^{(1)}_{p_{1}})
(L(-1)^{l_{1}}u_{(j)}^{(1)})_{j_{1}}
L(n^{(1)}_{1})\cdots L(n^{(1)}_{q_{1}})w_{(t)}^{1}\otimes}\nno\\
&&\quad\otimes L(-m^{(2)}_{1})
\cdots L(-m^{(2)}_{p_{2}})
(L(-1)^{l_{2}}u_{(j)}^{(2)})_{j_{2}}
L(n^{(2)}_{1})
\cdots L(n^{(2)}_{q_{2}})w_{(t)}^{2},\label{7-2},
\end{eqnarray*}
$m^{(1)}_{1}, \dots, m^{(1)}_{p_{1}}$, $m^{(2)}_{1}, \dots, m^{(2)}_{p_{2}}$,
$n^{(1)}_{1}, \dots, n^{(1)}_{p_{1}}$, $n^{(2)}_{1}, \dots, n^{(2)}_{p_{2}}
\in {\Bbb Z}_{+}$,
$l_{1}, l_{2}\in {\Bbb N}$, $j_{1}, j_{2}\in {\Bbb Q}$, $t=1, \dots, r$,
$j=1, \dots, k$.
Using the $L(-1)$-derivative property for generalized modules,
we see that elements of the form (\ref{7-2}) are spanned by
elements of the form (\ref{7-2}) with $l_{1}=l_{2}=0$. Now consider
the elements of the form (\ref{7-2}) with $l_{1}=l_{2}=0$
and of a fixed weight $s$. Since $W$ is lower-truncated, there are only
finitely many of them. This proves that the homogeneous subspaces of
$W$ are finite-dimensional. So $W$ is a $V$-module.

{\small \sc Department of Mathematics, Rutgers University,
New Brunswick, NJ 08903}

{\em E-mail address}: yzhuang@math.upenn.edu

\end{document}